\documentclass[12pt,a4paper,oneside]{amsart}

\usepackage{amsfonts, amsmath, amssymb, amsthm, amscd, hyperref}
\usepackage{anysize}
\usepackage{graphicx}
\usepackage[all]{xy}
\usepackage[utf8]{inputenc}

\newtheorem{theorem}{Theorem}[section]
\newtheorem*{theorem*}{Theorem}
\newtheorem{lemma}[theorem]{Lemma}

\newtheorem{corollary}[theorem]{Corollary}
\newtheorem{proposition}[theorem]{Proposition}

\theoremstyle{definition}
\newtheorem{definition}[theorem]{Definition}

\theoremstyle{remark}

\newtheorem*{remarks*}{Remarks}
\newtheorem{question}[theorem]{Question}

\numberwithin{equation}{section}

\DeclareMathOperator{\vol}{vol}

\DeclareMathOperator{\id}{id}

\begin{document}

\setcounter{page}{1}

\title{Lipschitz homotopies of mappings from 3-sphere to 2-sphere}

\author{Aleksandr Berdnikov}

\begin{abstract}
This work focuses on important step in quantitative topology: given homotopic mappings from $S^m$ to $S^n$ of Lipschitz constant $L$, build the (asymptotically) simplest homotopy between them (meaning having the least Lipschitz constant). The present paper resolves this problem for the first case where Hopf invariant plays a role: $m = 3$, $n = 2$, constructing a homotopy with Lipschitz constant $O(L)$.
\end{abstract}

\maketitle

\section*{Acknowledgments}

I am grateful to my advisor Larry Guth for introducing to me the interesting and rich topic of Lipschitz homotopies and for the guidance and help throughout the study; to Alexey Balitskiy for several corrections; and to Fedya Manin for developing a powerful theory from a few specific examples I considered here.

The research was supported by Larry Guth's Simons Investigator grant.

\setcounter{section}{-1}
\section{Introduction}
\label{sec:intro}

\subsection{Setting}

In this paper we consider Lipschitz mappings from $S^3$ to $S^2$ with standard metrics and homotopies between such mappings. Our main result is the following theorem. 

\begin{theorem}
\label{mainth}
Let $f_0,f_1: S^3\rightarrow S^2$ be mappings with Lipschitz constant $L$ that are homotopic to each other. Then there exists a homotopy $f_\bullet: S^3\times[0,1]\rightarrow S^2$ between them with Lipschitz constant $O(L)$.
\end{theorem}

Our theorem (and the methods introduced here) are an important step in the progress on the program proposed by Gromov in~\cite{gro} for understanding the interplay between hopotopical and metric complexities of mappings. One aspect of it is answering the following question:

\begin{question}
For a metric complex $X^m$ and homotopic mappings $S^n\to X^m$ with Lipschitz constant $L$, what minimal Lipschitz constant would be always sufficient for a homotopy (with such constant) to exist between them?
\end{question}

In the context of homotopies it is reasonable to separate the Lipschitz constant into spatial $|\partial /\partial x|$ and temporal $|\partial /\partial t|$ components. So we'll say that a homotopy that is $A$-Lipschitz in time variable and $B$-Lipschitz in space variable is, for shorthand, $A_t\times B_x$-Lip.

One of the first non-trivial considerations related to this question is attributed to Gromov in~\cite{gro}:

\begin{theorem}
\label{equi}
For $L$-Lipschitz mappings $S^n\to S^n$ one can always find $O(L_t\times L_x)$-Lip homotopies. 
\end{theorem}

The first step of untangling these questions is to reduce the mappings to a standard form. Since they ``look simple'' on the scale $1/L$, we can subdivide the domain into $1/L$-sized cells and bring a mapping to one of a few standard sample mappings on each cell. 

\begin{itemize}
\item For the case of $S^n\to S^n$ the procedure is especially straightforward. The mapping over the $n-1$-skeleton can be null-homotoped so that the Lipschitz constant keeps being $O(1)$ in a scaling where the cells are unit cubes. So the mapping is non-trivial only over now-separate $n$-cells, and bounded Lipschitz constant means that the mapping is close to some sample mapping --- a simplicial w.r.t. a subdivision that may be fine but is fixed --- and there are only finite number of those. 
\end{itemize}

This tweak is pretty much a more metric version of Brouwer's original approach where one boils down a mapping to a bunch of points in the domain, each carrying a mapping of degree $\pm 1$ around it. And a null-homotopy\footnote{We can restrict to just null-homotopies due to $\pi_n(X)$ being a group, and a homotopy $f\sim g$ being thus not much different from a null-homotopy of $f^{-1}g$} of such mapping is also just a more metric version of the Brouwer's argument, that canceling out the points of opposite signs, connecting them with disjoint curves in the homotopy domain. Originally that shows that the total degree is the only invariant on $\pi_n(S^n)$ thus proving it being $\mathbb{Z}$.

\begin{itemize}
\item In our case, to extend the mapping form $S^n$ to the domain of the homotopy $S^n\times[0,1]$, we want the canceling pairs of balls $B^n$ to become the opposite ends of embedded cylinders (instead of mere curves) that connect them in the bulk of the homotopy and cancel them out. In more applied terms we want to connect given sockets (the balls) with wires (the cylinders). The metric properties translate to the thickness of the wires and the size of the room their tangle occupies.\\
\item This problem (in a more general setting) was solved in~\cite{KB} by placing cylinders one at a time and showing that at each step the previous wires cannot block enough of possible tracks for the new one. Our approach provides an alternative solution, resembling a multi-stage election process: on $i$-th step groups of $2^i$ sockets interact to promote a neat output for the next scale to handle easily.  
\end{itemize}

Having done this construction, one gets a proof of the Theorem~\ref{equi}.

This approach uses the clear structure of maps $S^n\to S^n$, so that they can be dealt with in an explicit way by shuffling around the degrees in the standard little balls. In general we don't have this level of clarity about $\pi_m(S^n)$. But actually we don't need it either for a successful ``shuffling around of the degree''. The reason is, as long as the rest of the topological complexity constitutes just a {\itshape finite} group, it can tag along without occupying much space. That was one of the insights of~\cite{nullc}:

\begin{theorem}
The fact that $\pi_m(S^n)$ is finite for $m>n$ when $n$ is odd or $m<2n-1$ implies that in these cases there exist $O(L_t\times L_x)$-Lip homotopies for $L$-Lipschitz mappings $S^m\to S^n$. 
\end{theorem} 

We can build a required homotopy as follows. Subdivide $S^m\times [0,1]$, as always, into $1/L$-scale cubes. We build the map skeleton by skeleton. At the $n$-skeleton (the lowest interesting one) we need to figure out the mapping on each $n$-cube ${\bf C}^n_i$, i.e. the degree of the mapping ${\bf C}^n_i\to S^n$. Thus steps $n$ and $n+1$ are  mostly the same as in $S^n \to S^n$ case. But with each higher step we might first need to make some re-adjustments, otherwise the maps defined thus far may lack further extensions at all, when put together.

\begin{itemize}
\item For a cell ${\bf C}_i^k$ our mapping $\partial {\bf C}_i^k\to S^n$ might not extend inside ${\bf C}^k_i$, that is, represent a nontrivial element of $\pi_{k-1}(S^n)$. The hope is to tweak the responsible poor choices on the previous step, so that these obstructions vanish.\\
\item A-priori any tweaking might affect the Lipschitz constant that we care about. But if the whole set of possible obstructions, $\pi_{k-1}(S^n)$, is finite, then the set of ways to tweak a $(k-1)$-cell would be finite, as well as what can be happening on a $k$-cell in general, so we keep the bound ($O(1)$ on the scale of  a single cell).\\
\end{itemize}

\subsection{Formalization of the techniques}

The described approach promises more generality. When the geometry and topology get translated into a more algebraic language of cochains, it becomes easier to handle. Here are the parallels for a mapping $f: S^m\to X$.
\medskip
\begin{center}
\begin{tabular}{ |p{6cm}|p{6cm}| }
\hline
{\bf Geometry-wise} & {\bf Algebra-wise} \\ 
\hline
\medskip
$f|_{{\bf C}^k_i}$ on $k$-cells ${\bf C}_i^k\in S^m$, \newline the mapping to be extended \medskip & \medskip
$\phi: {\bf C}^{k+1}_i \mapsto [f|_{\partial {\bf C}^{k+1}_i}]\in \pi_k (X)$, \newline the cochain of obstructions \\ 
\hline
\medskip
$g_i\in \pi_k(X)$, the corrections \newline added to the $f|_{{\bf C}^k_i}$ \medskip & \medskip $\gamma : {\bf C}^{k}_i \mapsto [g_i]\in \pi_k (X)$, \newline the cochain of corrections \\
\hline
\medskip
Extendability (to a cell ${\bf C}_i^{k+1}$) \medskip &\medskip  $\phi + d\gamma = 0$ (on a cell ${\bf C}_i^{k+1}$)  \\
\hline
\medskip
Lipschitz constant remains \newline
bounded on the cell scale \medskip &\medskip sup-norm of $\gamma$ \newline is bounded \\
\hline
\end{tabular}
\end{center}
\medskip

The coefficients of cochains $\phi$ and $\gamma$ --- the homotopy groups --- reflect the geometric complexity of the problem. But their finite (torsion) aspects, as we've seen, don't matter much for the bounds in the end. That suggests to focus just on {\itshape free} homotopy invariants. They are much simpler and can be computed through pullbacks of differential forms. 
\begin{itemize}
\item The prime example of such invariant would be the usual degree of a mapping $f: S^n\to S^n$. It can be computed via $\int f^* d\vol_{S^n}$.\\
\item A more interesting example would be the integral Hopf invariant of a mapping $f: S^{2n-1}\to S^n$. It can be computed via $\int \alpha \wedge \beta$ where $\alpha =f^* d\vol_{S^n}$, and $\beta$ is an anti-derivative of $\alpha$ ($\alpha=d\beta$).\\
\item Other invariants like these, that can be computed (in a simply connected case) via $\int \bigwedge \eta_i$ where $\eta_i$ are obtained from pull-backs of volume forms by products and taking anti-derivatives. 
\end{itemize}

This intuition has solidified in the form of the {\itshape shadowing principle} of~\cite{nullhlin}. It shows, how to return to the geometry once we are done with the algebra, that is, how to build back an actual homotopy (or any mapping) given only the (cochains of) pullbacks of some key forms. When one boils a space down to the essential structure of differential forms, one gets Sullivan minimal models $\mathcal{M}^* (X)$ --- simplest DGAs that acceptably mimic DGA $\Omega^*(X)$ of forms on $X$.

\begin{theorem}
(Shadowing principle, sketch summary). Given a formal mapping $\phi : \mathcal{M}^* (Y)\to \Omega^*(X)$, one can recover a genuine mapping $f: X\to Y$, so that $f^*: \Omega^*(Y)\to \Omega^*(X)$ is close to $\phi$.
\end{theorem}

\begin{itemize}
\item The formal mapping tells us how some important forms (like $f^* d\vol$) we'd like to be arranged in the domain $X$ and the DGA structure ensures that it isn't done in a contradictory way. 
\item From this perspective, shadowing principle is simply a powerful generalization of what was happening in the case $S^n\to S^n$. In that problem the algebra $\mathcal{M}^*(S^n)$ (the ``formal $S^n$'') is generated by a single element $d\vol$, the generator of $H^n(S^n)$, if we look only up to degree $2n-2$.
\item Thus, to build a formal mapping $\mathcal{M}^{< 2n-1}(S^n) \to \Omega^*(X)$ means to fix the closed form $\alpha = f^* d\vol$.
\item It is sufficient for $\alpha$ to be a (universally bounded) cellular cochain on a fine enough subdivision of $X$, where the values $\alpha({\bf C}^n)= \int_{\bf C^n} f^* d\vol$ of the cochain on cells ${\bf C}^n$ are the degrees of the resulting mapping on these cells.
\item Summarizing this example, to extend (in our case from $X\times \{0,1\}$ to $X\times [0,1]$) an $L$-Lipschitz mapping to a sphere $S^n$ it is enough to extend $f^* d\vol$ by a closed form (given that the dimension of the space is at most $2n-1$). If sup-norm of the extension is bounded by $\sim L^n$, the Lipschitz constant of the extended mapping would be bounded by $\sim L$.     
\end{itemize}

This important tool seals away a lot of technicalities that tend to repeat in this topic. It allows the arguments to easily address the substance without getting lost in details.

\begin{itemize}
\item In the case of building a homotopy $h$ between a pair of mappings $f,g:X\to Y$, one may pick a general algebraic formula for a formal homotopy $\eta : \mathcal{M}^* (Y)\to \Omega^*(X\times [0,1])$ between the formal mappings $\phi,\gamma : \mathcal{M}^* (Y)\to \Omega^*(X)$ representing $f,g$ --- and get back an actual homotopy $h$.
\item In simpler cases one only have to deal with forms that represent cohomology (like $\alpha = f^* d \vol$). They can be extend as $\eta(d\vol)=(1-t)\alpha$, so that they dissipate throughout the homotopy. 
\item But the degrees cannot just vanish on their own, they are a conserved quantity --- or, algebraically, $d\eta(d\vol)$ should  be closed, since $d\vol$ itself was closed, but $d(1-t)\alpha=-\alpha\neq 0$. So we need to find where to move the degrees so that they cancel out --- find for $\alpha$ an anti-derivative $\beta$ so that $d\beta=\alpha$. Then move the degrees $\alpha$ along $\beta$ --- that is, define $\eta(d\vol):=(1-t)\alpha - dt\wedge \beta$, so that $d\eta(d\vol)=0$. 
\item Finding such $\beta$'s (small enough anti-derivatives) is, at this point, an ever-present step in the argument. 
Not all exact forms have anti-derivatives that are just as small. But if we fix a finite complex $X$, then for any exact form $\alpha\in \Omega*^(X)$ we can find $\beta$ such that $d\beta =\alpha$ and $||\beta||\leqslant C_{\bf IP}$ where $C_{ \bf IP}(X)$ is a finite {\it isoperimetric constant} of $X$.
\end{itemize}

This approach\footnote{At that point the shadowing principle wasn't introduced yet; rather, the principle is the generalisation of the techniques and experience from these papers.} was used in~\cite{nullc} and~\cite{nullhom} to get more general estimates:

\begin{theorem}
Mappings to a simply-connected space, rationally equivalent to a product of Eilenberg-MacLane spaces, posess $\sim L_t\times L_x$-Lip homotopies~\cite{nullc}.
\end{theorem}

\begin{theorem}
Mappings to a simply-connected homogeneous space (eg. $S^n$) have $\sim L^2_t \times L_x$-Lip null-homotopies~\cite{nullhom}.
\end{theorem}

These were improved in \cite{nullhlin}:

\begin{theorem}
Mappings $X\to Y$ to any simply-connected space have $\sim 1_t \times L^2_x$-Lip {\it null-homotopies} and $\sim 1_t \times L^p_x$-Lip {\it general homotopies} for some $p(X,Y)$. If $Y$ is symmetric, null-homotopies can (for any $\varepsilon>0$) be made $\sim L^{1+\varepsilon}$ Lipschitz~\cite{nullhlin}.
\end{theorem}

\subsection*{This paper and beyond}

The results listed in the end of the previous section mostly employed the following reasoning. Take the mappings to be homotoped, consider their formal versions, interpolate between them uniformly, and then try to patch up any introduced inconsistencies. Such straightforward approach runs into the following issues. Some mappings, like Hopf mapping, are significantly better at storing large degrees in a tight space than more simple ones, because of the non-local interactions exhibited in Hopf invariant: the total is more than the sum of the parts for the Hopf invariant. So we might have great trouble dealing with these secretly dense mappings if we manage them uniformly throughout the time, bit by bit. 

\begin{itemize}
\item For example, the Hopf mapping $H$ is 3-dimensional, but scales like $L^4$. A way to see it is to compose it with an $\sim L$-Lipschitz self-map $F$ of $S^2$, such that $\deg(F)= L^2$. One way to compute Hopf invariant of a mapping $f$ is to extend $\alpha f^* d\vol$ to a ball $B^4$ and integrate $\int_{B^4}\alpha^2$. $F$ multiplies $\alpha$ by $L^2$, so the degree becomes $L^4$.
\item So, if we approach such a package armed only with unit Hopf mappings that we can put on our cells, we lose this compression and are bound to spend $L$-fold more time (giving $L^2_t$ results for certain spheres).
\item Specifically, consider a mapping $S^3\to S^2$ composed of an $L^4$ degree mapping on one half of $S^3$ and its reflection of degree $-L^4$ on the other half. If we were to null-homotope that mapping using only Hopf mappings on 3-cells, it will take a lot of time. Indeed, $L^4$ of Hopf invariant would have to pass across the dividing equator. But it is 2-dimensional and so has only $\sim L^2$ cells and lets $L^2\times L$ unit Hopf invariant maps to pass.
\end{itemize}

The solution is to shuffle these mapping without unpacking them. That is exploited in the $L^{1+\varepsilon}$ result above, where smaller-scale mappings are iteratively used to build mappings for the next scale.

In this paper we resolve the issue of such compressed degrees with a trick of similar recursive flavor. At each step we group the current cells of the domain $S^3$ into bigger cells of double the size. New cells may aggregate more complexity of a mapping inside, than the initial basic units. So we comb new cells and morph the mapping to resemble one of the few basic mappings. At which point an iteration of this process can be applied again, and so forth until the victory.

Such approach allows to handle the ``over-dense'' regions of the mappings directly. From a certain moment such piece is viewed as a single indivisible unit, so it doesn't matter what's in the box --- we open it no more. This optimizes the required time down to $\sim L$, matching the lower bound. This paper is considering a specific case $S^3\to S^2$ (a bit too explicitly and by hand), but the same arguments relayed in the obstruction language give the same construction of linear time for all cases of $S^m\to S^n$ and more.

In~\cite{scal} this process of coarsening of a mapping (and hence, building $\sim L$-Lipschitz homotopies) was generalized to a broad class of spaces $Y$, dubbed {\itshape scalabe}. Name comes from the defining feature of having self-maps $Y\to Y$ that scale $Y$ well, in the sense that they act on $H^n(Y,\mathbb{R})$ by $L^n$ while having Lipschitz constant $\sim L$. If we can scale our domain in this sense, then we can build nice coarser basic mappings, pass to more and more coarse cells and get $L_t\times L_x$ homotopies just as before. 

Why bother with scaling maps, though? For a start, scalability is easier to check due to its another defining feature. It appears when one considers scaling maps of higher and higher degrees. Since they have linear Lipschitz bounds, their pullbacks of forms (after normalization by $L^{-\deg}$) are bounded. And in the limit these forms give a presentation of cohomology $H^\bullet (Y)\to \Omega_\flat (Y)$ (morally, because $f*d$ is scaled out). And vica versa, such presentation is a guide for shadowing principle to construct scaling mappings --- and hence, efficient homotopies: 

\begin{theorem}
For a (formal simply connected finite) complex $Y$ the existence of $\sim L_t\times L_x$ nullhomotopies in $[X,Y]$ for any (finite) $X$ is equivalent to the existence of a presentation $H^\bullet (Y)\to \Omega_\flat (Y)$ of DGA of its cohomology (sending each class to its representative)~\cite{scal}.
\end{theorem} 

Thus, one can easily check scalability for, say, connected summs of $\mathbb{C}P^2$'s or $\mathbb{H}P^2$'s or $\mathbb{O}P^2$'s and find which admit linear homotopies\footnote{For example, $\#_n \mathbb{O}P^2$ is 
scalable iff $2n\leqslant {16 \choose 8}  = 12870$.} 

Currently we pursue $L\log^{p(X)}(L)$-Lipschitz homotopies for some of non-scalable formal spaces. They do have self-maps that scale homology, but without metric bounds, so we construct them rather manually. The construction ended up an optimisation of the recursive construction from~\cite{nullhlin}.

\section{Overview}

The general structure of the proof is as follows. 

First we replace the source space $(S^3,*)$ with a cube $([0,1]^3,\partial [0,1]^3)$ that is $\sim 1$-biLipschitz embedded in the sphere $S^3$. 

Our method is based on the notion of ``cubical mappings'' that are simple on some scale $\mathcal{L}$. We start by homotoping $f_i$ to mappings that are cubical on a scale, corresponding to the Lipschitz constant of the mappings. The aim now is to repeatedly pass to coarser grids by merging cells together, while maintaining the simple cellular structure of the mapping. Thus, the key lemma for the proof is the following statement.

\begin{lemma}
\label{main}
Let $n$ be a positive integer and $f:([0,2n\mathcal{L}]^3,\partial [0,2n\mathcal{L}]^3)\rightarrow (S^2,*)$ be a $\sim 1$-Lipschitz mapping that is cubical on the grid of size $\mathcal{L}$. Then there is a homotopy of length ${\sim}\mathcal{L}$ connecting $f$ through ${\sim}1$-Lipschitz mappings to a mapping that is cubical on the size $2\mathcal{L}$. 
\end{lemma}

This lemma allows to conclude the proof of Theorem~\ref{mainth}. We begin with mappings that are cubical on a unit scale. Applying repetitively Lemma~\ref{main} we end up with mappings cubical on the scale of the entire domain. The simple structure of cubical mappings (knowing that the boundary $\partial [0,2n\mathcal{L}]^3$ is sent to the point) guaranties then that two such mappings are bounded homotopic. The total length of the combined iterated homotopies is then bounded by $\sim \mathcal{L}$ and all the intermediate mappings have Lipschitz constant ${\sim}1$.

The process can be roughly illustrated by the following picture in a simpler case of mappings $S^1\to S^1$.

\medskip

\includegraphics[scale=0.2]{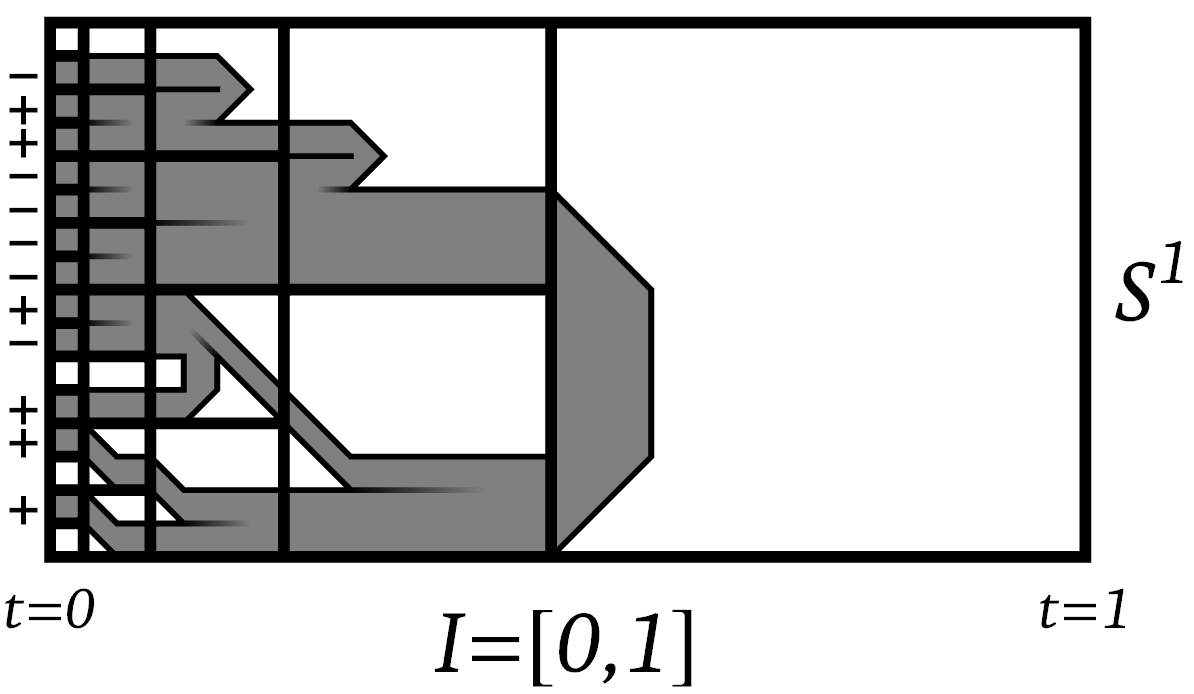}

The picture represents the domain of the homotopy $[0,1]\times S^1$. At $t=0$ one has $S^1$ divided in $\sim \mathcal{L}$ segments, each carrying one of 3 standard mappings: of degree $+1$, $-1$ or 0. On each step pairs of neighboring segments merge into greater ones and bundle their standard mappings into greater standard ones (with the degrees bounded by $2^k$ on step $k$). The mappings are represented by the (bundles of) gray wires connecting the original elementary mappings of the opposing degree (for more details see Definition~\ref{wire}). The whole homotopy is messy, but each highlighted square is simple because it deals only with 2 incoming mappings and 1 outgoing, so all square can be dealt with in a uniform manner.  

In our case, the proof of Lemma~\ref{main} is where the most of the technicalities reside. It requires to fix some burdensome terminology to even sketch the argument, so for details we refer the reader to Section~\ref{proof}. 

\section{Hopf Links}

In this section we introduce some elementary blocks used in the construction and learn to operate with them. We want to think about a mapping of a manifold $M^m$ into $S^n$ in terms of its fibers in $M^m$. That is, we use a Lipschitz analog of the construction of Pontrjagin and Thom that builds a map $f: M^m \to S^n$ from a framed submanifold $W^{m-n}\subset M^m$. The submanifold $W$ is meant to be the preimage $f^{-1}(s)$ of a regular value in $s\in S^n$, and the framing of $W$ lets trivialize its tubular neighborhood as $D^n\times W$. To build the map $f$ one defines it on the neighborhood $D^n\times W$ by projecting onto the disk factor and collapsing its boundary, and sending the rest of $M$ to the basepoint. 

To handle the Lipschitz constant in this construction, we work with the tubular neighborhood $D^n \times W$ itself, rather than a single fiber $W$. Here is a more precise definition.

\begin{definition}
\label{wire}
A {\bf wire} $l$ in $M^3$ is a ${\sim}1$-biLipschitz embedding $l:D^2\times (W,\partial W)\rightarrow (M^3,\partial M^3)$ where $W$ is a connected 1-manifold (with or without boundary), and $D^2$ is a unit disk. The corresponding {\bf wire mapping} given by $l$ is the composition
$$M^3\supset \text{Im}(l)\stackrel{l^{-1}}{\rightarrow} D^2\times W\stackrel{pr}{\rightarrow} D^2\stackrel{/\partial D^2}{\rightarrow} S^2.$$
\end{definition}

We would like to pack wires into cables. 

\begin{definition}
\label{defcable}
A {\bf cable} $c$ (of size $d$ and degree $N$) of constant cross-section is a collection of $N$ wires $l_i$ given by the compositions
$$D^2\times (W,\partial W)\stackrel{d_i\times \id}{\longrightarrow}B^2\times (W,\partial W)\stackrel{c}{\rightarrow} (M^3,\partial M^3),$$
where $d_i$ are disjoint translates inside of a disk $B^2$ (of diameter $\leqslant d$) of ${\sim}1$-biLipchitz embedding of $D^2$ into a square $[0,1]^2$, and $c$ is a ${\sim}1$-biLipchitz embedding into $M^3$. The {\bf cable mapping} (of constant cross-section) given by $c$ is the mapping given by all of its wires altogether.
\end{definition}

For the purposes of this section these definitions will be enough, but we will have to expand it in Section~\ref{proof} to deal with the case when the cross-section changes along the cable. Until then all cables are assumed to be of constant cross-section.

We defined mappings on the wires so that they send the boundary of the mapping to the base point. So we may extend the mapping to the whole $M^3$ by sending the rest of $M^3$ to the base point. Also, the homotopies of such mappings are meant rel boundary, meaning the restriction to the boundary is constant. 

We can now lock the cables into ``Hopf links'' to define mappings with some non-trivial Hopf invariant:

\begin{definition}
Denote by ${\bf a}|{\bf b}$ the Hopf link of cables of degrees $a$ and $b$ (and size $d\sim \sqrt{a}+\sqrt{b}$). Its cable mapping is the Whitehead product of mappings $S^2\rightarrow S^2$ of degrees $a$ and $b$.
\end{definition}

This is one of the fundamental building blocks of our construction: it carries Hopf invariant $H(a|b)=2ab$ inside a cube with edge length ${\sim}d\sim \sqrt{a}+\sqrt{b}$. Given an even Hopf invariant $2n$ we can efficiently represent it by a link $a|b+c|1$. Here $a=\lfloor \sqrt{n}\rfloor$, $b=\lfloor n/a \rfloor$ and $c=n-ab$, so that $a,b,c \sim \sqrt{n}$ and $H(a|b+c|1)=2n$. By {``+''} in this formula we mean a mapping defined on two neighboring cubes of size ${\sim}\sqrt{n}$, on one cube by $a|b$, on the other --- by $c|1$. This gives a way to carry some Hopf invariant; but what we need to do, is not just allocate fixed portions of Hopf invariant here and there, but freely operate on them, splitting and merging during the homotopy. In this section we create machinery that will achieve this.

During this section we will use the notation $A\backsimeq B$ to indicate that a sum $A$ of $a|b$'s is homotopic to another sum $B$ via a ${\sim}1$-Lipschitz homotopy of length that is $\sim$ to the linear sizes of $A$ and $B$. Note that as long as the number of summands is bounded, we don't need to specify how exactly summands are located and oriented. Indeed, any two possible arrangements can be homotoped into each other in time $\sim$ to the size of the mapping: one just needs to rotate and move to new location each of the $a|b$'s, which takes linear time. The specific details of this procedure are straightforward. Also, we don't have to worry about the specifics of the deg-$a$ and deg-$b$ mappings used to define the construction, since any two such mappings are homotopic to each other in linear time by Proposition 1.6 of~\cite{larry} (or much more general Theorem B in \cite{nullhom}).

We start our construction with the following geometric observation and will proceed in a more algebraic fashion.

\begin{proposition}
\label{1geom}
There exists a following homotopy of cable mappings: $$a_1|b+a_2|b\backsimeq (a_1+a_2)|b.$$
\end{proposition}

\begin{proof}
The homotopy is illustrated by the following figures:

\includegraphics[scale=0.15]{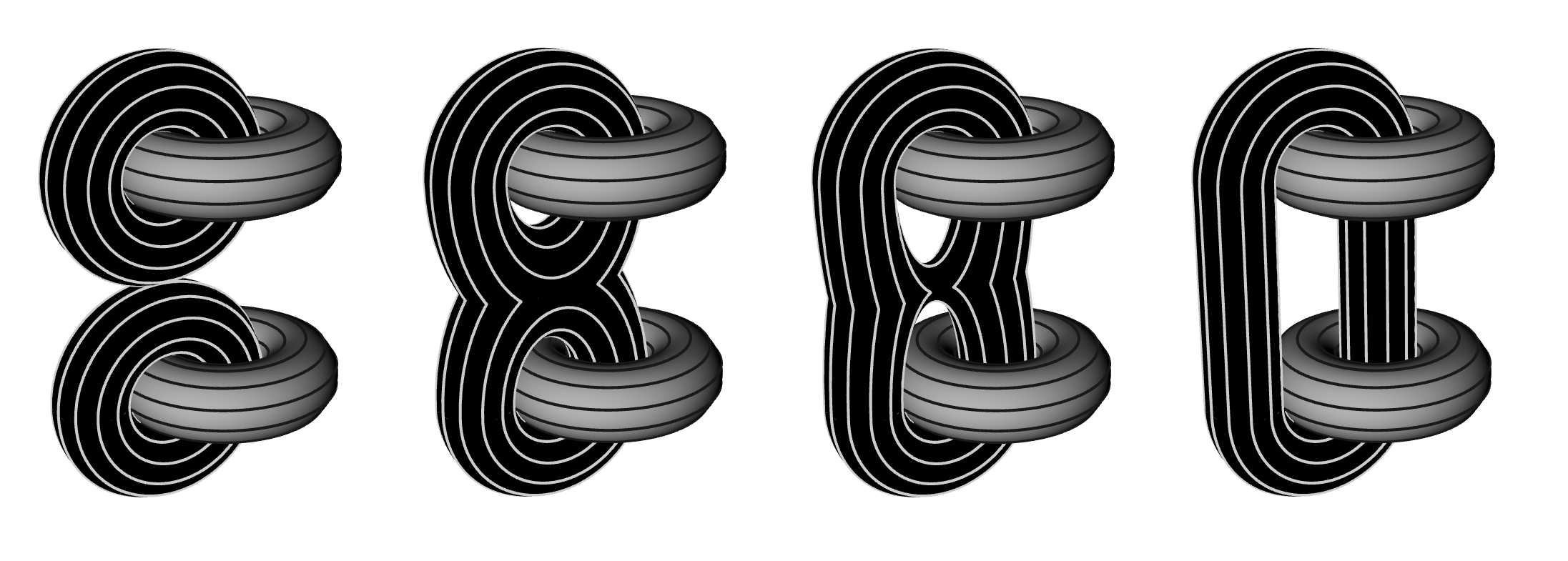}

Here the light tori denote cables of degree $a_i$ and the dark tori --- cables of degree $b$. We start with two $a_i|b$'s located on top of each other so that the $b$-cables are symmetric about the horizontal plane $z=0$. That is, the part of the mapping $B_0(x,y,z)$ described by the $b$ cables depends on $x,y$ and only $|z|$ (instead of $z$). Let the centers of links $a_i|b$ have $|z|=z_0$. Then at time $t$ redefine the mapping $B$ in the region $|z|<z_0$ by the formula for the homotopy 
$$B_t(x,y,z)=B_0(x,y,z_0+(z_0-t)(|z|/z_0-1)).$$ 
This homotopy has the same Lipschitz constant in time direction as in the spatial direction and at time $t=z_0$ it merges two $b$-cable mappings into a single one. So we ended up with $(a_1+a_2)|b$ in time $z_0$.
\end{proof}

\begin{corollary}
\label{2bal}
There exists a following homotopy of mappings: $$2a|b\backsimeq a|b +a|b \backsimeq a|2b.$$
\end{corollary}

Now we use these moves to arrange the links so that they use the domain efficiently. Indeed, the order of the size of a link $(a|b)$ depends only on the greater of $a$ and $b$, so unless they are of the same order, we use the space not efficiently. To remedy that, we ``balance'' our links:

\begin{definition}
A link $a|b+c|1$ is \textbf{balanced} if $2b\geqslant a\geqslant b,c$ and its size is ${\sim}\sqrt{a}$. 
\end{definition}

\begin{proposition}
\label{crude}
For a mapping $a|b+c|1$ such that $a\geqslant b$, $a\gtrsim c$, there exists a homotopy between mappings $a|b+c|1\backsimeq \alpha|\beta + \gamma| 1$ so that $ab+c=\alpha\beta + \gamma$ and the latter link is balanced.
\end{proposition}

\begin{proof} 
First we can reduce to the case $a>c$ by applying the homotopy 
$$a|b+c|1\ \ \backsimeq\ \ a|(b+1)-a|1+c|1\ \ \backsimeq\ \ a|(b+1)+(c-a)|1.$$
until $c$ becomes smaller than $a$. Now all we have to do is to essentially apply repeatedly Corollary~\ref{2bal} until $b$ grows to become comparable to $a$. More specifically, let the link that we have on the $n$-th stage be $a_n|b_n+c_n|1$. If $a_n$ is odd, do
$$
\begin{array}{c}
a_n|b_n+c_n|1\ \ \backsimeq \ \ (a_n-1)|b_n+1|b_n+c_n|1\ \ \backsimeq
\vspace{0.3cm}\\
\backsimeq \ \ (a_n-1)|b_n+(c_n+b_n)|1\ \ =\ \ a'_n|b'_n+c'_n|1.
\end{array}$$
Now, when $a'_n$ is even, employ Corollary~\ref{2bal}:
$$
\begin{array}{c}
a'_n|b'_n+c'_n|1\ \ =\ \ 2a'_{n+1}|b'_n+c'_n|1\ \ \backsimeq
\vspace{0.3cm}\\
\backsimeq \ \  a'_{n+1}|2b'_n+c'_n|1 \ \ =\ \ a'_{n+1}|b'_{n+1}+c'_{n+1}|1.
\end{array}$$
If after that $c'_{n+1}\geqslant a'_{n+1}$, reduce $c$ again as in the beginning of the proof. This action needs to be taken only 3 times at most during each step since by induction we had $a_n>c_n$ and $a'_{n+1}>(a_n-1)/2-2$. After we achieve that $c$ is less than $a$, those are our new $a_{n+1}|b_{n+1}+c_{n+1}|1$.

Now, iterate this procedure. Each time it doubles $b$ and halfs $a$ up to an additive constant $C$, so by repeating the procedure we can achieve $a_N\leqslant 2b_N+3/2C$ while $a_N\geqslant b_N$. Use $c|1$ several times (${\sim}C$) as before to transfer the degree from $a_N$ to $b_N$ in order to get to $a_N\leqslant 2b_N$. The resulting link $a_{N+1}|b_{N+1}+c_{N+1}|1$ is the goal link $\alpha|\beta + \gamma| 1$.

The length of homotopy of the $n$-th iteration is ${\sim}\sqrt{a_n}\leqslant \sqrt{a/2^n}$, so the total length of the homotopy assembled from these iterations is at most
$${\sim}\sum_{n=1}^\infty \sqrt{a/2^n}=\frac{1}{1-1/\sqrt{2}}\sqrt{a}\sim \sqrt{a}.$$ 
\end{proof}

This proposition allowed one to balance a link in a crude way: we changed $a$ and $b$ by factors of 2 --- in giant steps --- and therefore only managed for $a$ and $b$ to become of the same order of magnitude. But what we would like, is to be able to rebalance links in any way we may need:

\begin{proposition} 
\label{prec}
For any two balanced links $a_i|b_i+c_i|1$ of equal Hopf invariant there exists a homotopy of mappings $a_1|b_1+c_1|1\backsimeq a_2|b_2+c_2|1$.
\end{proposition}

\begin{proof} $$a_1|b_1+c_1|1 \ \ \backsimeq \ \ (a_1|b_1+c_1|1) +\big(-(a_2|b_2+c_2|1)+ (a_2|b_2+c_2|1)\big)= $$
$$=\big( (a_1|b_1+c_1|1) -(a_2|b_2+c_2|1)\big) + (a_2|b_2+c_2|1)\ \ \backsimeq \ \ a_2|b_2+c_2|1,$$
where the first homotopy is trivial and the second one is a consequence of the following proposition.
\end{proof}

\begin{proposition}
\label{cancel}
For any two balanced links $a_i|b_i+c_i|1$ of equal Hopf invariant there exists a homotopy of mappings $a_1|b_1+c_1|1-a_2|b_2-c_2|1\backsimeq 0|0$.
\end{proposition}

\begin{proof} 
First, balance $-a_2|b_2-c_2$ and $a_1|b_1+c$ using Lemma~\ref{crude} (w.l.o.g. $a_i\geqslant b_i$ after that). Then, merge $c_i$'s into a single $(c_1-c_2)|1=:c|1$ (w.l.o.g. $c>0$). 

Now we will iteratively lower all the numbers until they become less than a (later chosen) constant $B$ that doesn't depend on the size of the links we start from. There exist only a finite amount of choices for $0<a_i,b_i,c<B$, so we may pick a bounded null-homotopy for each choice (such that the total Hopf invariant is 0). Once the numbers get below $B$, we apply that null-homotopy, and due to a finite pool it is uniformly bounded. And to get to that point we repeat the following  procedure:
\begin{enumerate}
\item Using homotopies $\pm a_i|b_i + c|1\backsimeq \pm a_i|(b_i\mp 1) + (c+a_i)|1$ arrange so that $b_i$ are even.
\item Now lessen the quantities involved by partially canceling them:
$$a_1|b_1-b_2|a_2\ \ \backsimeq\ \ (a_1-\frac{b_2}{2})|b_1+\frac{b_2}{2}|b_1-b_2|\frac{b_1}{2}-b_2|(a_2-\frac{b_1}{2})\ \ \backsimeq $$
$$\backsimeq \ \ (a_1-\frac{b_2}{2})|b_1-b_2|(a_2-\frac{b_1}{2})\ \ =:\ \ a'_1|b'_1-a'_2|b'_2.$$
\item Balance $-a'_2|b'_2$ and $a'_1|b'_1+c$, swap $a$ and $b$ if necessary to get $a_i\geqslant b_i$.
\end{enumerate}

Note that at each step the quantity $a_2b_2$ decreases at least by a constant factor during the step (2). Indeed, as long as $a_2>B$, all $a_i$ and $b_i$ are of the same order of magnitude, i.e. their ratios are bounded by, say, 3; indeed, we can achieve this with some value of bound $B$, so we set $B$ to that value. Once we know $b_1\geqslant a_2/3$, we get 
$$\frac{a'_2b'_2}{a_2b_2}=\frac{a_2-b_1/2}{a_2}\leqslant \frac{a_2-a_2/6}{a_2}=5/6.$$
Each cycle consists of $<10$ homotopies whose length is ${\sim}\sqrt{d_n}$ where $d_n$ is the highest degree involved in the step $n$, $$d_n\sim(a_2)_n\leqslant \sqrt{(a_2b_2)_n}\lesssim \sqrt{(a_2b_2)_0}(5/6)^n<(a_2)_0(5/6)^n,$$
so the total length is less than
$$\sim \sum_{n=0}^{\infty}10\sqrt{d_n}\lesssim \sum_{n=0}^{\infty}\sqrt{(a_2)_0(5/6)^n}=\frac{1}{1-\sqrt{5/6}}\sqrt{a_2}\sim \sqrt{a}.$$
\end{proof}

Now that we can do arbitrary rebalancing, we can freely add our links as in the next lemma.

\begin{lemma}
\label{ad}
For any two balanced links $a_i|b_i+c_i|1$, with the first having higher Hopf invariant, there exist homotopies between mappings
$$a_1|b_1+c_1|1+a_2|b_2+c_2|1\ \ \backsimeq \ \ d_1|e_1+f_1|1,$$
$$a_1|b_1+c_1|1-(a_2|b_2+c_2|1)\ \ \backsimeq \ \ d_2|e_2+f_2|1$$
for any balanced links $d_i|e_i+f_i|1$ that match the Hopf invariant.
\end{lemma}

\begin{proof} 
For the addition split $d_1|e_1+f_1|1$ into two links with the same $d_1$ that have the same Hopf invariants as the $a|b+c|1$-summands. Use Lemma~\ref{prec} to rebalance them into $a_i|b_i+c_i|1$. For subtraction analogously do
$$a_1|b_1+c_1|1\ \ \backsimeq \ \ d_2|e_2+f_2|1+(a_2|b_2+c_2|1)$$ 
and then cancel the last term with its negative.
\end{proof}

\section{Interlocked and twisted Hopf links}

The links studied in the previous section are simplest tools to manage the Hopf invariant of a mapping. However, the mappings we deal with usually have their Hopf invariant stored in a more convoluted way. In this section we introduce generalized versions of links $a|b$. They will serve as a transients between the complexity of mappings emerging during the homotopy and the simplicity of $a|b$ links.

\subsection{Interlocked links}
\ \newline
 
First version is the ``interlocked link'' $(a_i)\begin{smallmatrix}|\\n \end{smallmatrix}(b_i)$:

\begin{definition}
Given sequences $(a_i)$ and $(b_i)$, $1\leqslant i\leqslant n$ and $\max (a_i)=a$, $\max(b_i)=b$, define the \textbf{ interlocked link} $(a_i)\begin{smallmatrix}|\\n \end{smallmatrix}(b_i)$ as cables $A$ and $B$ constructed as follows (see the picture below):
\begin{enumerate}
\item Cross-section of cable $A$ consists of $2n$ stripes $A_i$ of shape $1\times a$.
\item Stripe $A_{2i-1}$ is empty (to leave space for interlocking $B_i$) and $A_{2i}$ is a row of $a_i$ mappings of degree 1 (that is, cross-sections of corresponding wires). Same goes for $B$.
\item Cables $A$ and $B$ are rings that intersect across a box ${[0,a]\times [0,b] \times [0,2n]}$ so that the stripe $A_{2i}$ passes through $[0,a]\times [0,b]\times {[2i-1,2i]}$ along the second coordinate and $B_{2i}$ passes through $[0,a]\times [0,b]\times {[2(n-i),2(n-i)+1]}$ along the first coordinate.
\item The rings $A$ and $B$ close up at different sides of their intersection, so that for any $k$ the cables formed by $a_{k<i\leqslant n}$ and $b_{n-k+1\leqslant i\leqslant n}$ form a Hopf link and cables formed by $a_{1\leqslant i\leqslant k}$ and $b_{1\leqslant i< n-k}$ aren't linked at all.
\end{enumerate}   
\end{definition}

This definition is illustrated by the following interlocked link $(8,8,8,8)\begin{smallmatrix}|\\4 \end{smallmatrix}(8,8,8,8)$:

\includegraphics[scale=0.25]{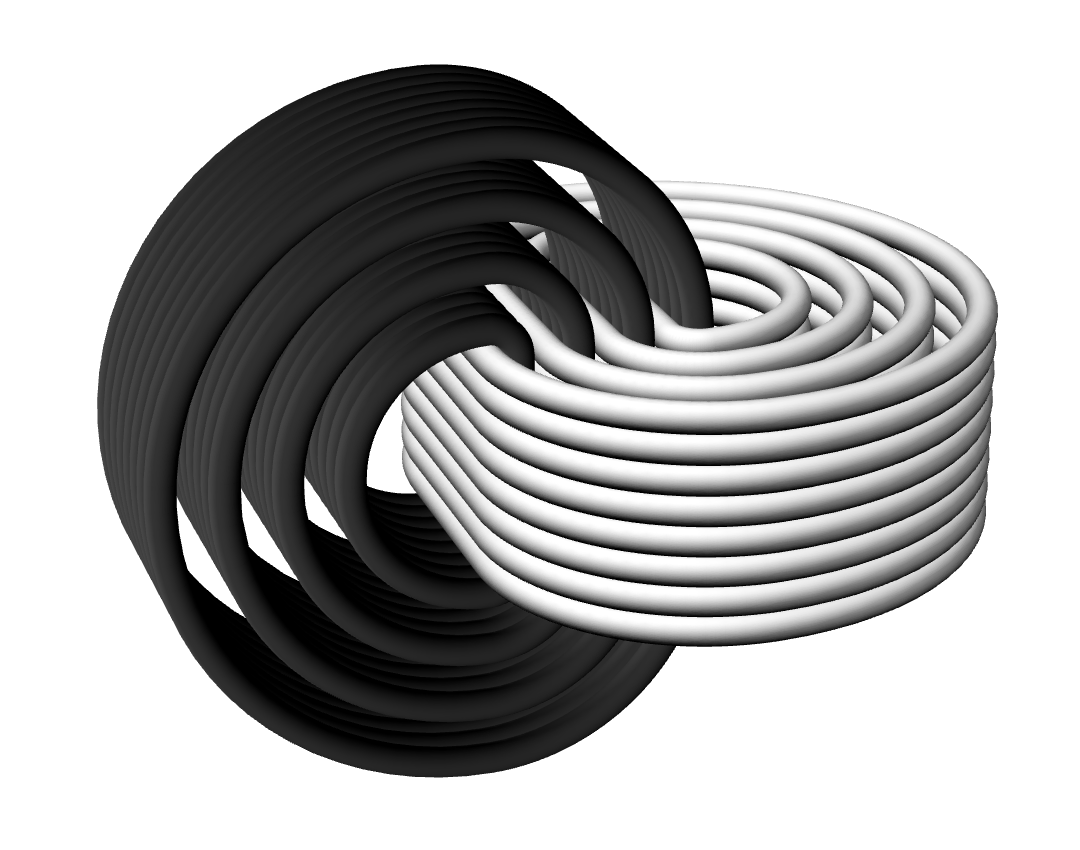}

It is clear that the size of the mapping $(a_i)\begin{smallmatrix}|\\n \end{smallmatrix}(b_i)$ is ${\sim}\max (a_i,b_i,n)$ (where by the size we mean the linear dimensions of the smallest box encompassing the image of the cable mapping).

As the purpose of such mappings is to convert more complicated ones into the standard links $a|b$, the main lemma regarding the interlocked links is the following.

\begin{lemma}
\label{inter-bal}
For any link $(a_i)\begin{smallmatrix}|\\n \end{smallmatrix} (b_i)$ and a balanced link $a|b+c|1$ of the same Hopf invariant there is a homotopy between mappings
$$(a_i)\begin{smallmatrix}|\\n \end{smallmatrix} (b_i)\backsimeq a|b+c|1.$$
\end{lemma}

\begin{proof} First add to the interlocked link $(a|b+c|1)-(a|b+c|1)$, so that we only need to cancel the second term with the interlocked link. The statement of the lemma now follows from iteratively applying the next proposition. 

\begin{proposition}
\label{interprop}
For any link $(a_i)\begin{smallmatrix}|\\n \end{smallmatrix} (b_i)$ and a balanced link $a|b+c|1$ of the same Hopf invariant there is a homotopy of mappings
$$(a_i)\begin{smallmatrix}|\\n \end{smallmatrix} (b_i)-a|b-c|1\ \ \backsimeq\ \ \sum_{k=1}^8 A_k,$$
where each $A_k$ is a mapping of the form
$$\Bigl(\Bigl[\frac{a_i}{2}\Bigr]\Bigl)\begin{smallmatrix}|\\ [\frac{n}{2}] \end{smallmatrix} \Bigl(\Bigl[\frac{b_i}{2}\Bigr]\Bigr)-a'|b'-c'|1,$$
that satisfies the assumptions made for the initial mapping; $[x]$ denotes here rounding to either closest integer.
\end{proposition}

The proof is postponed until Lemma~\ref{inter-bal} is taken care of. The proposition allows to iteratively split interlocked links of the lemma into 8 smaller ones with all the quantities roughly halved each time. Locating the 8 new links (together with their standard link partners) in cells of a ${2\times 2 \times 2}$ subdivision of the initial cell we get the same setting as in the beginning but with quantities reduced by half (more precisely, if $a_i,b_i,n,\sqrt{a},\sqrt{b},\sqrt{c}$ are bounded by $2^k$, then the corresponding quantities of the descendant links are bounded by $2^{k-1}$). Applying Proposition~\ref{interprop} until all quantities become $<10$ we get to a point where all the mappings are from a finite list of null-homotopic mappings (at most one for each choice of $n,a_i,b_i,a,b,c<10$). Choosing a Lipschitz null-homotopy for each mapping from the list we conclude by applying it to the last stage of iteration. That step is bounded uniformly and others are bounded as before by a convergent geometric series (with factor $1/2$), so the total length is bounded by ${\sim}$~the~length of the first iteration which is bounded by Proposition~\ref{interprop}.
\end{proof}

{\bf Proof} of Proposition~\ref{interprop}. The algorithm is best described alongside an example. Here we display subsequent stages of it for the case of $(8,8,8,8)\begin{smallmatrix}|\\4 \end{smallmatrix}(8,8,8,8)$ that we used before, omitting the standard link $a|b+c|1$. Remark that in the general case all the quantities $a_i,b_i,c_i$ are of the order of the linear size of the mapping. In the example we have $a_i=b_i=8$ and $n=4$, but these specific values will not be used in the argument. 

At the first step of the algorithm we deflate the cables in half so that both sequences $(a_i)$ and $(b_i)$ are split roughly in the middle by a gap of size ${\sim}n$:  

\includegraphics[scale=0.3]{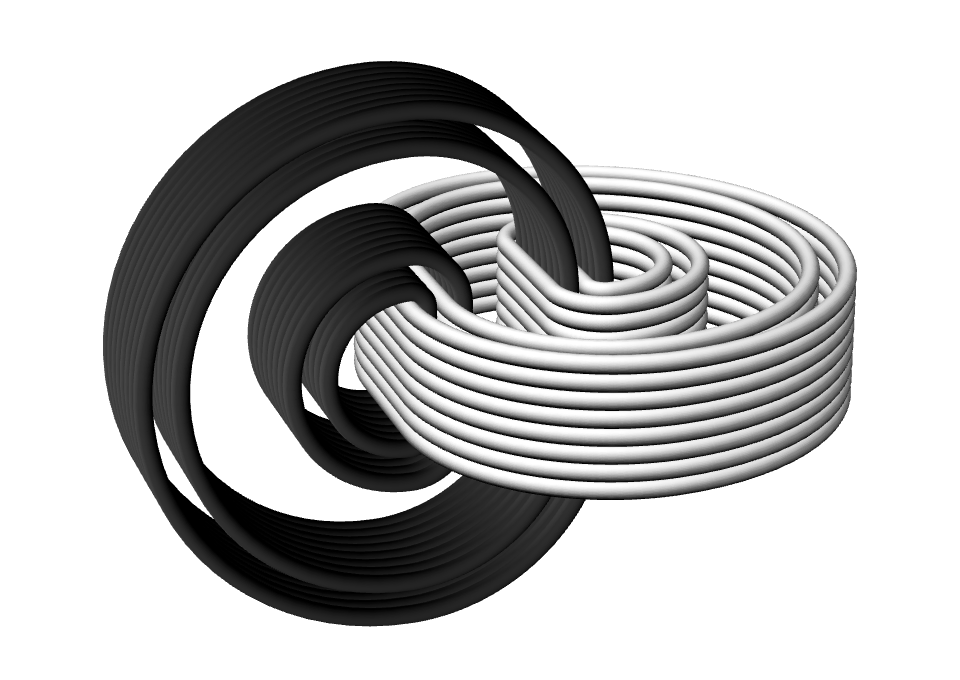}

That is done so that we can independently rotate the interlocked links $(a_{1\leqslant i\leqslant \lceil \frac{n}{2}\rceil})\begin{smallmatrix}|\\ \lceil \frac{n}{2}\rceil\end{smallmatrix} (b_{\lfloor \frac{n}{2}\rfloor <i \leqslant n})$ (the small black sub-cable with the large white one) and $(a_{\lceil \frac{n}{2}\rceil<i\leqslant n})\begin{smallmatrix}|\\ \lfloor \frac{n}{2}\rfloor\end{smallmatrix} (b_{1\leqslant i \leqslant \lfloor \frac{n}{2}\rfloor})$ (the large black one with the small white one). So we rotate them as a whole around the axes of larger cables so that the smaller ones get to the opposite side:

\includegraphics[scale=0.3]{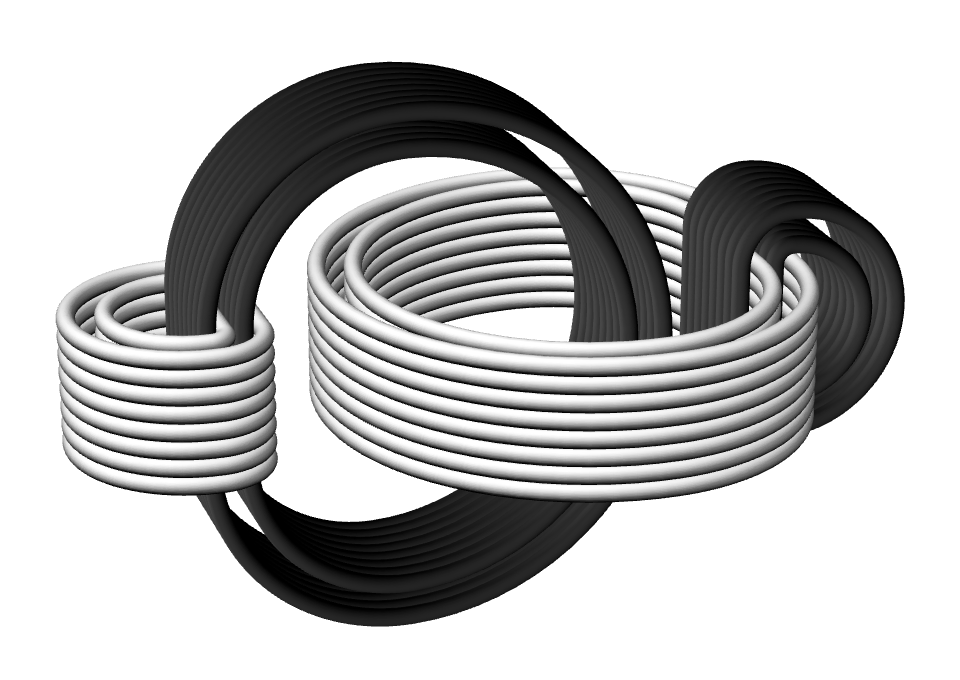}

We do it so that now the larger cables interact with each other and with smaller cables independently. Because of that a procedure similar to the one used in Proposition~\ref{1geom} can be used to split the larger cables into two smaller ones: one linked in a standard way with a part of other large cable, and one interlocked with its smaller cable:

\includegraphics[scale=0.3]{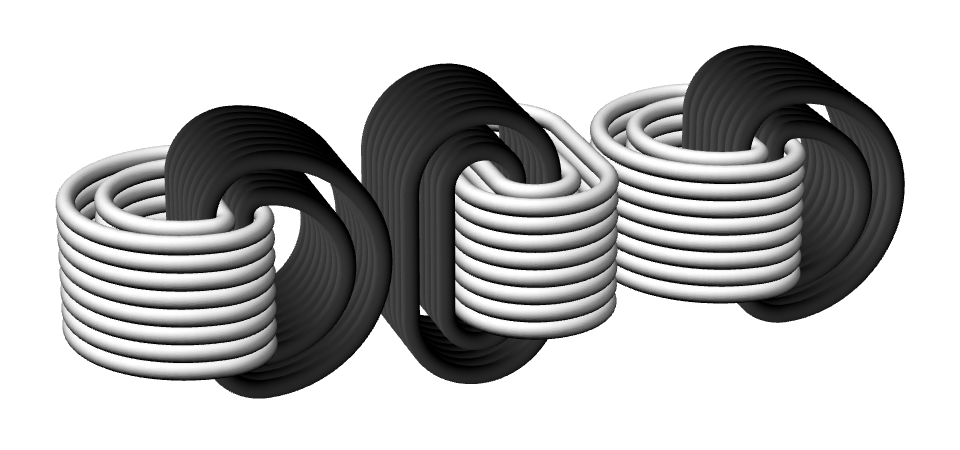}

In other words, we get two interlocked links on the sides and a standard one in the middle. The standard one merges next with (or, cancels a part of) the standard link $-(a|b+c|1)$ from the statement of the proposition, so we are left with two interlocked links. They have $a_i$ and $b_i$ of the same magnitude, but roughly twice as small $n$. Then $a_i$ and $b_i$ can be divided in half as well using homotopies as in Proposition~\ref{1geom} (construction of which for the interlocked case is practically the same):

\includegraphics[scale=0.2]{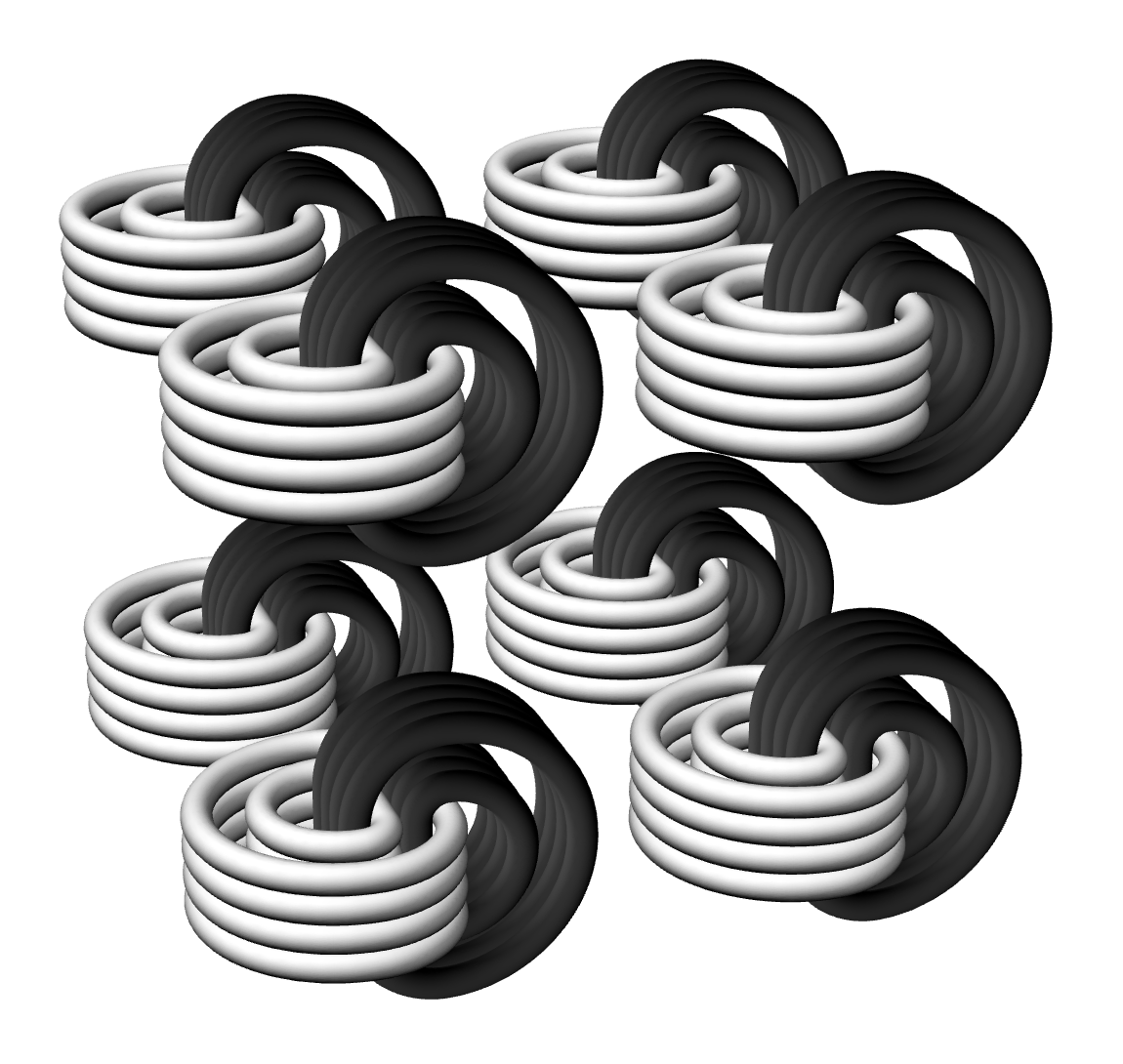}

Once we end up with 8 smaller interlocked links, we split the remains of the standard link $-(a|b+c|1)$ into eight balanced links that match Hopf invariants of the corresponding interlocked links (total Hopf degree was 0 so those 8 pieces precisely exhaust the remains of $-(a|b+c|1)$).

\begin{flushright}
$\square$
\end{flushright}

\subsection{Twisted links}

\ \newline

The job of the interlocked links is to be an intermediate step in breaking up the ``twisted links" --- another sort of building block that will occur in the course of the homotopy.

\begin{definition}
A \textbf{twisted link} $\underline{a}$ is a mapping that one gets from a standard link $a|0$ by altering it by a half of a Dehn twist:
\newline
\includegraphics[scale=0.2]{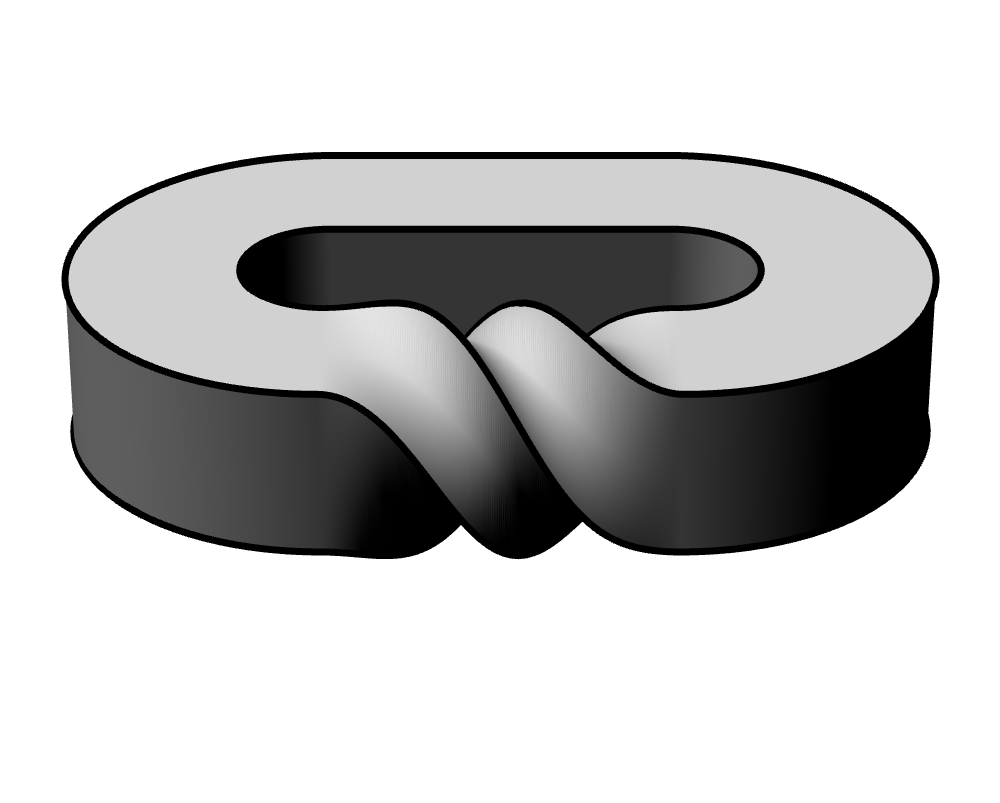}
\end{definition}

{\bf Explanation.} Let the link $a|0$ be such that all the cross-sections $C_a$ consist of $|a|$ identical $\deg$-1 mappings $f:(D^2,\partial D^2)\rightarrow (S^2,pt)$, whose centers $c_i$ are located in a centrally symmetric fashion. Let $[0,L]\times C_a$ be a segment of the link of length $L$ that is of the same order as the linear sizes of the cross-section $C_a$. Consider a ${\sim}1$-biLipschitz map $T$ of $[0,L]\times C_a$ into its neighborhood that on a cross-section $\{t\}\times C_a$ is a rotation by $\pi\frac{t}{L}$. To get a twisted link $\underline{a}$ redefine the mapping on $[0,L]\times C_a$ so that on the section $\{t\}\times C_a$ the mappings $f$ are centered not at $c_i$ but at $T(c_i)$. Note that the mappings maintain their orientation and are {\itshape not} rotating together as a whole $C_a$ (otherwise on the end $t=L$ we wouldn't return to the same mapping because there was only half-turn rotation).

Another piece of new notation comes from the necessity to deal with possibly odd Hopf invariants while all the constructions up to this point carried an even Hopf invariant. Denote by $\overline{1}$ a fixed ${\sim}1$-Lipschitz mapping of size ${\sim}1$ with Hopf invariant 1 realized by a wire loop with a 1-turn Dehn twist. Also we write $\overline{(1)}$ to denote a link $\overline{1}$ coming with coefficient either of $\{-1,0,1\}$ to avoid considering different parities of Hopf invariant separately. Fix a homotopy $\overline{1}+\overline{1}\backsimeq 1|1$. It allows to modify all our algorithms regarding balanced links $a|b+c|1$ carrying arbitrary {\itshape even} Hopf invariant so that to generalize to links $a|b+c|1+\overline{(1)}$ that can now carry {\itshape any} Hopf invariant. The modification constitutes transferring occasional excesses of $\overline{1}$-s into $c|1$-terms of a standard link by means of the homotopy $\overline{1}+\overline{1}\backsimeq 1|1$ that we fixed.

Now we can tie our new objects to the previous ones:

\begin{proposition}
\label{twist-inter}
For a twisted link $\underline{\alpha}$ there is a homotopy of mappings
$$\underline{\alpha}\backsimeq (a_i)\begin{smallmatrix}|\\n \end{smallmatrix} (b_i)+a|b+c|1+\overline{(1)}$$
where all the numbers on the right hand side are bounded by ${\sim}$the size of $\underline{\alpha}$. 
\end{proposition}

\begin{proof}
First of all, group the centers $c_i$ into circular families $C_i={c_{i,k}}$. For example, say that $c_{i,k}$ are distinguished by having the distance from $c_{i,k}$ to the center $O$ of the cross-section of $\underline{\alpha}$ being in the range $(i-1,i]$. By slightly homotoping the mapping we can assume that the distance is actually of the same value $i$ across each family $C_i$ and the centers $c_i$ are located equidistantly across the ring.

Now the twisted region splits into similar blocks. First, there are tubular shells, each of them containing one family $C_i$. Then each such shell is cut into rings, each containing not the full turn of the family $C_i$ but rather carrying each $c_i$ to the next one. Each such ring in turn consists of segments like in the left picture:
	
\includegraphics[scale=0.2]{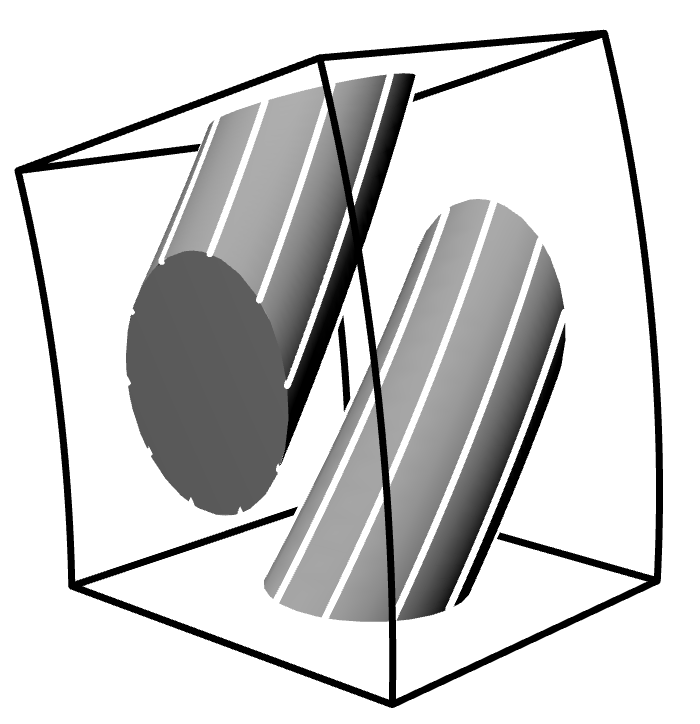}
\includegraphics[scale=0.2]{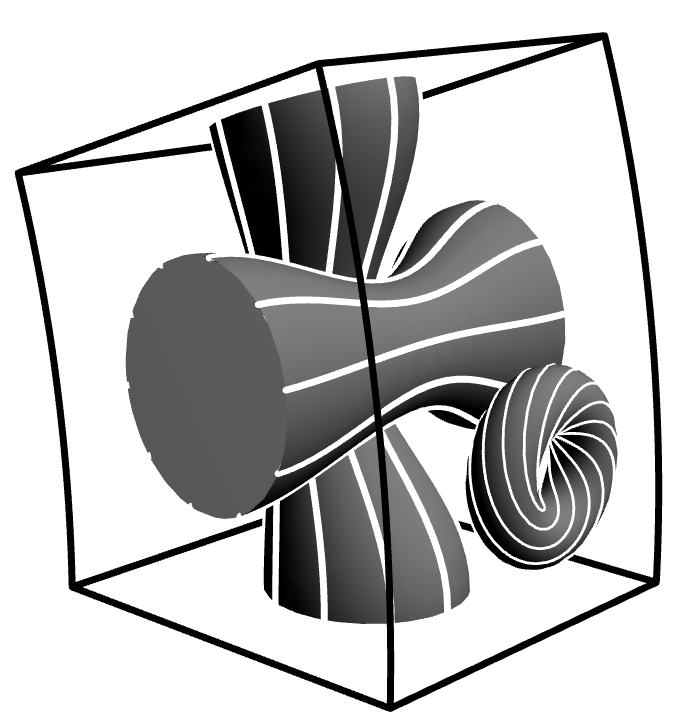}

It is bounded-homotopic (relative boundary) to the mapping on the right, where the torus is the $\overline{1}$-mapping needed to even out the (relative) Hopf invariant. It doesn't really matter that the difference is exactly Hopf invariant 1, we only use that it is the same finite number for all the segments, so we don't actually calculate it. By applying such a homotopy to all the segments we transform the mapping into several parts. One is just now untwisted link $\alpha|0$, its part in the right picture is horizontal. The second part is a bunch of circular wires located in the bulk of $\alpha|0$, in the picture such wire is vertical. And finally, there is plenty of $\overline{1}$'s. The latter are moved out of the link and merged into a single balanced link $a|b+c|1+\overline{(1)}$ in a binary fashion, similar to that of Lemma~\ref{inter-bal}, but in reverse.

Now the cross-section in the (previously) twisted region looks (after minor tweaks) like the following picture on the left:

\includegraphics[scale=0.2]{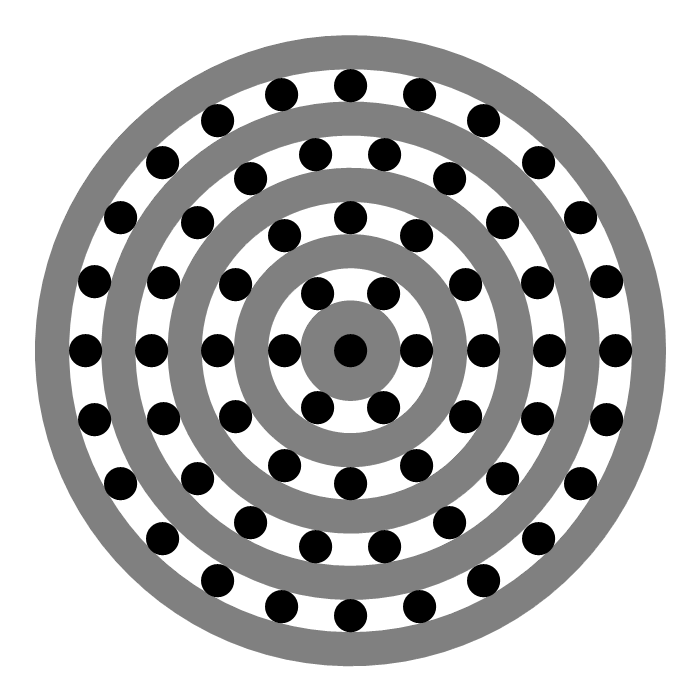}
\includegraphics[scale=0.2]{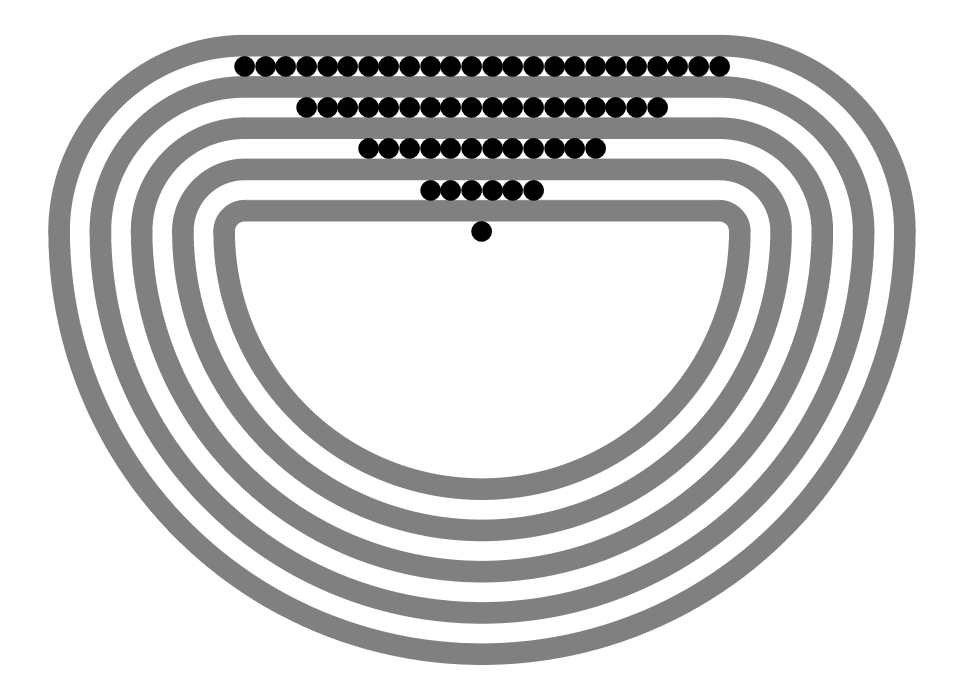}

Here black circles are cross section of wires of $\alpha|0$ (horisontal in the previous illustration) and gray rings are circular wires coming from the vertical pieces in the previous illustration. Homotope all the cross-sections of $\alpha|0$ according to the homotopy transforming the left picture to the right collecting all black disks to one side. After that the wires of $\alpha$ are interlocked with the gray circular wires.
\end{proof}

\begin{corollary}
\label{twist-bal}
For a twisted link $\underline{\alpha}$ there is a homotopy of mappings
$$\underline{\alpha}\backsimeq a|b+c|1+\overline{(1)}$$
where all the numbers on the right hand side are bounded by the ${\sim}$linear size of $\underline{\alpha}$. 
\end{corollary}

{\bf Proof} is a combination of Proposition~\ref{twist-inter} and Lemma~\ref{inter-bal}.

\section{Definitions}
\label{defs}

In order to infuse our arguments with some rigor, we have to refine our definitions to be more explicit and general.

\subsection{Definition of a general cable mapping}

\ \newline

Fix a ${\sim}1$-Lipschitz degree-1 mapping $w:(D^2,\partial D^2)\rightarrow (S^2,*)$. Recall that we call a ${\sim}1$-biLipschitz embedding $l:D^2\times (W,\partial W)\rightarrow (M^3,\partial M^3)$ a {\bf wire} in $M^3$ where $W$ is a 1-manifold (with or without boundary). Any such wire defines a mapping 
$$wl^{-1}:M^3\supset \text{Im}(l)\rightarrow S^2,$$ which we call a {\bf wire mapping} given by $l$. In the definition~\ref{defcable} of a cable we packed these wires together to form a mapping of higher degree. One limitation of that version of cable definition is that the cross-sections of the cable were made all the same throughout the cable, so that relative positions of wires are fixed. We would like to generalize the definition so that it encapsulates the possibility of wires varying their relative arrangement as they go along the cable. 

To maintain some structure on the arrangement of the wires in the new more flexible setting, we define the following intermediate object. An $n$-{\bf stripe} (or just {\bf stripe}) of length $a$ is a collection of $n$ wires put in a box $[0,n]\times[0,1]\times[0,a]$, with $i$-th wire being mapped by $c\times \id_{[0,a]}+(i,0,0)$ where $c:D^2\rightarrow [0,1]^2$ is a standard embedding.

Now we define an abstract cable.
In what follows $K$ is a fixed constant (such that $K\geqslant 2$) and it is used to thicken all the middle part of the cable (except for the very ends) to give enough room for stripes in the cable to move around. 
So, the {\bf abstract cable} of size $n$ and length $a$ is a collection of $\leqslant n$ disjoint stripe embeddings 
$$S_i: ([0,n_i]\times[0,1])\times[0,a]\rightarrow [0,Kn]^2\times [0,a],$$ 
$$\text{such that }S_i(x,\alpha)=(x+v_i(\alpha),\alpha),$$
for $v_i:[0,a]\rightarrow \mathbb{R}^2$ being ${\sim}1$-Lipschitz paths, such that the stripes are in a strict order, if to compare lexicographically coordinate-wise. More precise, we require that if $i<j$ then for all $\alpha\in [0,a]$ we have (with the superscript indicating the coordinate component)
$$v_i^2(\alpha)+1\leqslant v_j^2(\alpha)\text{, or }$$
$$v_i^1(\alpha)\leqslant v_j^1(\alpha)\text{ and }v_i^2(\alpha)\leqslant v_j^2(\alpha).$$
We also require that the ends $([0,n_i]\times[0,1])\times\{0,a\}$ of the stripes map into smaller faces $[0,n]^2\times \{0,a\}$.

A {\bf cable} is a ${\sim}1$-biLipschitz embedding $f$ of an abstract cable into some $M^3$. More accurately, $f$ is an embedding of the space $[0,Kn]^2\times [0,a]$ of the abstract cable. Then the compositions $f\circ S_i$ define embeddings of corresponding stripes, and therefore of the wires, into $M^3$. They in turn provide a mapping $\text{Im}(f)\rightarrow S^2$ as usual. This mapping given by all these wires is what we call ``the mapping given by this cable'', when extended by a constant mapping the rest of the domain.

The inequalities in the definition of an abstract cable mean that the stripes $S_i$ preserve their order in the bulk of a cable, meaning that they are first stacked in rows (along the first coordinate, just as wires in each stripe) and then those rows are stacked in a column (along the second coordinate). The purpose of this structure is, on one hand, to make the cross-section of a cable malleable enough, and on the other hand --- to keep stripes (and hence wires) from tangling inside the cable. Managing arrangements of the stripes is further simplified by the following argument.

\begin{proposition}
\label{comb}
Let $(S_i)$ and $(S_i')$ be two abstract cables of size $n$ of length $a\geqslant n$, such that $(S_i)|_{\{0,a\}}=(S_i')|_{\{0,a\}}$. Then there exists a ${\sim}1-$Lipschitz homotopy $(S_i)_t$ through cables with $(S_i)_0=(S_i)$ and $(S_i)_n=(S_i')$, and the homotopy is constant on the ends of the cables.
\end{proposition}  

There is a caveat: the intermediate cables may have the velocities $|\partial v_i/\partial \alpha|$ greater (by a constant factor) than the constant fixed for other cables. The reason is, we may somewhat increase the velocities during the homotopy, so if the initial cable had a maximal one, we will overshoot the usual bound. It is not an issue since these faster cables occur only temporarily during the homotopy and the final cable is a normal one.  

\begin{proof}
By a 2-Lipschitz reparametrization of the segment $[0,a]$ we may homotope both cables to be constant (and therefore equal) on $[0,a/4]$ and $[3a/4,a]$. 

Now perform the following homotopy. On the segment $[a/4,3a/4]$ it linearly drags the stripes $S_i$ and $S_i'$ along the second coordinate to make it equal to $i$. Extend this homotopy to the rest of the segment $[0,a]$ by linear interpolation to a constant homotopy. No overlapping of stripes (say, $S_i$ and $S_j$, $i<j$) occurs since either they were in the same row and hence disjoint by the first coordinate alone (since it doesn't change); or they were in different rows, but then both at the start and at the end of the homotopy we have $v^2_i+1\leqslant v^2_j$, so this inequality holds for all times, hence $S_i$ and $S_j$ are disjoint by the second coordinate in this case.

At this stage stripes $S_i$ and $S_i'$ have matching second coordinates on the whole $[0,a]$, matching first coordinates outside of $[a/4, 3a/4]$, and the second coordinate separates stripes in both families on $[a/4, 3a/4]$. Therefore the linear homotopy from $v_i$ to $v_i'$ will not overlap any stripes and hence concludes the homotopy we build.

Notice that at the first (reparametrization) step the velocities $|\partial v_i/\partial \alpha|$ are multiplied by 2, in the second step they gain up to $4K$ and the third step doesn't increase these velocities.
\end{proof}

Proposition~\ref{comb} can be viewed as a parametric/relative version of the following simpler proposition.

\begin{proposition}
\label{reshape}
Given a pair of 0-length cables $(S_i)^{start}$ and $(S_i)^{end}$ given by stripes
$$S_i^\bullet:[0,n_i]\times[0,1]\rightarrow [0,n]^2,$$
there exists an abstract cable $S_i$ over $[0,T]$ with $T{\sim}n$ such that 
$$S_i|_{\{0\}}=S_i^{start}\text{ \ and \ \ }S_i|_{\{T\}}=S_i^{end}.$$
\end{proposition}      

\begin{proof}
The construction, if $[0,T]$ is thought of here as time, is identical to the homotopy built in Proposition~\ref{comb} for a middle point ${a/2\in[0,a]}$ if first cable was given there by $(S_i)^{start}$, and the second --- by $(S_i)^{end}$. Since Proposition~\ref{comb} works independently for all $\alpha\in [a/4,3a/4]$, it is fine that in applying the proof of Proposition~\ref{comb} we defined the cables only for $\alpha= a/2$.
\end{proof}

\subsection{Definition of cubical mapping}

\label{cubical}

We start by defining a mapping cubical on scale 1. First, by (a Lipschitz analog of) simplicial approximation theorem there is a finite collection $S_{cell}'$ of mappings $[0,1]^3\rightarrow S^2$ sending the 1-skeleton $([0,1]^3)^{(1)}$ to a point, so that for some $c>0$ any $c$-Lipschitz mapping 
$$([0,n]^3,([0,1]^3)^{(1)})\rightarrow (S^2,pt)$$
is ${\sim}1$-Lipschitz homotopic (rel $([0,1]^3)^{(1)}$) to a mapping that on each unit cube of the grid restricts to a mapping from $S_{cell}'$, and also, each mapping from $S_{cell}'$ on each face of the cube has $|\deg|\leqslant 1$.

Once the pool of mappings of unit cubes is bounded, more specific features for $S_{cell}'$ will be secured. Namely, for any possible restriction $f^{face}$ to a face of a unit cube of a mapping from $S_{cell}'$ fix a (rel $\partial$) ${\sim}1$-Lipschitz homotopy $f_t^{face}$ of that mapping to a cross-section of a wire mapping. Then for any mapping $f'\in S_{cell}'$ fix a Lipschitz extension $f'_t$ of homotopies $f_t^{face}$ on the faces, so that the end mapping of the homotopy is given by several wires going from one face of the cube to another avoiding sub-cube $[1/4,3/4]^3$, and inside that sub-cube there is a balanced link that adjusts the Hopf invariant. The collection $S_{cell}$ of these end mappings encompass what we define to be \textbf{cubical} on scale 1 mappings. By the construction we get that the initial $c$-Lipschitz mapping is ${\sim}1-$Lipschitz homotopic to a mapping $f_1$ that is cubical on scale 1.

Now given a mapping $f_1:[0,n\mathcal{L}]^3\rightarrow S^2$ (for some positive integer $n$) that is cubical on scale 1, we would like to build a cubical version $f_\mathcal{L}$ of $f_1$ that has similar simple structure at the scale $\mathcal{L}$, replacing wires by cables. The next picture illustrates the kind of look a cubical mapping can have on one $\mathcal{L}$-cell. It shows the simplicity of organized structure we aim at. However for illustrative purposes the picture deviates from minor and arbitrary specifics of the construction in the paper.

\includegraphics[scale=0.33]{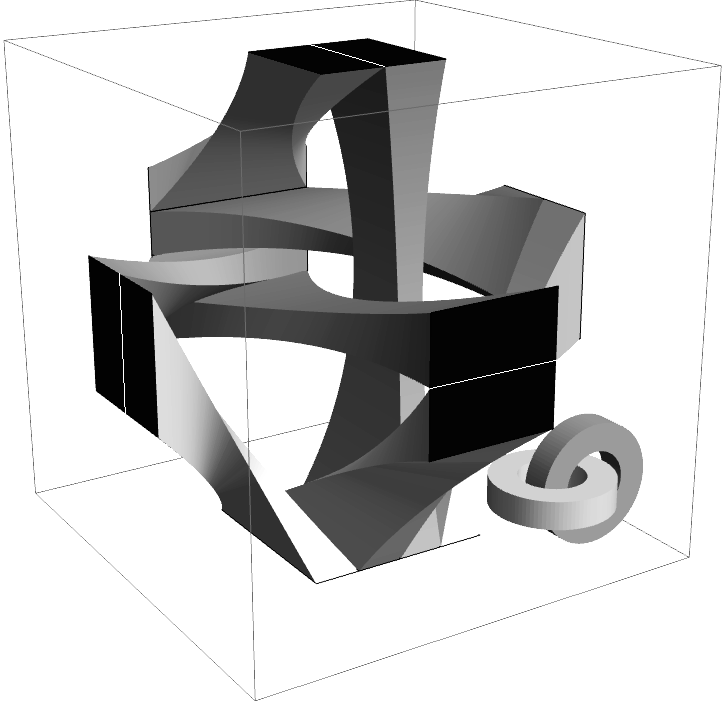}

We cannot quite literally pick the wires of $S_{cell}$ as guidelines for cables since we may have to pass several cables through a single face. Because of this and some other issues we introduce a template that will be used later to guide cables. We take a unit cube $[0,1]^3$, and for each coordinate plane fix order and directions of coordinates. The \textbf{template} is a fixed collection of ${\sim}1$-biLipschitz cable template embeddings 
$$C_i:([0,1]^2\times[0,1])\cup([0,K]^2\times[1/4,3/4])\rightarrow [0,1]^3$$ 
such that
\begin{itemize}
\item Cable templates $C_i$ are in one-to-one correspondence to the ordered pairs of faces of the cube $[0,1]^3$,
\item thus, to each cable $C_i$ two faces are assigned, a positive $P_i$ and a negative $N_i$.
\item For each template $C_i$ its face $[0,1]^2\times \{0\}$ is mapped to $N_i$ and the face $[0,1]^2\times \{1\}$ --- to $P_i$. 
\item In the coordinates on $P_i$ and $N_i$ each mapping in the previous item is either the identity mapping or identity with the second coordinate flipped; the choice is made so that $C_i$ end up preserving the orientation near each face.
\item For each face the templates meeting that face are ordered.
\item The templates don't intersect the region $[1/4,3/4]^3$.
\item The templates satisfy a non-intersection condition.
\end{itemize}

The last condition need some explanation but roughly means that the cables built by these templates cannot intersect. It would be easy to achieve by just picking the templates themselves to be non-intersecting, but that is not an option since each face of the cube may belong to several templates, so they necessarily intersect there. Hence, to satisfy the ``non-intersection condition'' one has to rely on the specifics of how the cables are put through templates. Therefore the proof of existence of the cable is given in Proposition~\ref{promise} that is postponed until the construction of the actual cables is determined.

Now, presuming existence of a template, we explain how to build from it a cubical on scale $\mathcal{L}$ version $f_\mathcal{L}$ of a mapping $f_1$ cubical on scale 1. We start with 2-skeleton of $\mathcal{L}$-grid. On each face $[0,L]^2$ carrying degree $d$ we populate the first $|d|$ unit sub-cells with a fixed wire cross-section $c$ and the rest --- with a constant mapping. The sub-cells are compared coordinate-wise (like in definition of an abstract cable, second coordinate trumps the first) with respect to coordinates of templates ends. We choose either $P_i$ or $N_i$ so to make the resulting map of the same sign as $f_\mathcal{L}$ and presume this chart on each face from now on.

Given $f_\mathcal{L}$ on the faces, we extend it inside the cell by a few cables factoring through the (rescaled) templates. Assign a degree $d_i\geqslant 0$ to each cable $C_i$ of the template so that on each face each cable contributes a non-negative degree and the total matches the degree $d$ of $f_\mathcal{L}$. Then, if the first cable going through the face is $C_i$, assign to it the first $d_i$ of unit sub-cells of the face, then $d_j$ sub-cells after that --- to the second cable $C_j$, etc. Extend the cables in a constant way through templates until the middle half of templates and fill in the rest by applying Proposition~\ref{reshape}.

More precisely, in order to use~\ref{reshape} we have to group the unit sub-cells on the faces (that form the ends of the cable, $S_i^{start}$ and $S_i^{end}$) in $\lesssim K\mathcal{L}$ stripes. We do it in a greedy way, declaring two wires to belong to the same stripe whenever it is possible: that is, if on both ends they have the same second coordinate (and hence are located in the same row and potentially the same stripe). Then in one face each row of cells is broken by such conditions (coming from the other face) into at most 2 pieces, since the breaks it generates are spaced by $\mathcal{L}$ cells apart from each other. So, there are no more than $2\mathcal{L}\leqslant K\mathcal{L}$ stripes in the cable, as required.

Now it is a good time to quickly fulfill our promise on non-intersecting conditions of the templates.

\begin{proposition}
\label{promise}
There is a template such that the cables (built by it as above) never intersect.
\end{proposition}

\begin{proof}
First, to define the mappings on the first quarter of the cable length, $[0,1/4]$, extend the templates from each face. Along this extension make them quickly shrink $100K$-fold towards the center, but also linearly spread apart (along the first coordinate) in their respective order so that the faces $[0,K]^2\times \{1/4\}$ map into the same plane but don't intersect. Since the index of the cable (that a wire belongs to) is monotone with respect to the first coordinate, no intersection of wires of different cables will occur by this point. Treat the last quarter, $[3/4,1]$, symmetrically. Secondly, for the inner half of the templates, $[1/4,3/4]$, just connect its ends (that were made disjoint in the previous step) by disjoint fillings that avoid $[1/4,3/4]^3$.
\end{proof}

Returning back to $f_\mathcal{L}$, we reserved its inner part $[\mathcal{L}/4,3\mathcal{L}/4]^3$ for a balanced link that would equalize Hopf invariants of $f_1$ and $f_\mathcal{L}$ on each $[0,\mathcal{L}]^3$ cell. That does not quite make sense yet, since these mappings disagree on the boundary of the cell. To resolve the ambiguity, we fix a ${\sim}1$-Lipschitz over time $\mathcal{L}$ homotopy $h_\mathcal{L}$ on each face $[0,\mathcal{L}]^2$, so that $h_\mathcal{L}$ connects (rel boundary) $f_1$ to $f_{\mathcal{L}}$. Now we know which balanced link should be put inside the 3-cell so that the homotopy $h_\mathcal{L}$ from $f_1$ to $f_\mathcal{L}$ extends from $\partial [0,\mathcal{L}]^3$ to the interior. Such link is defined by having Hopf invariant that is opposite to the invariant of the rest of the the mapping. That is, the one that is given on the sides of $\partial ([0,\mathcal{L}]^3\times[0,T])$ by $h_\mathcal{L}$ and on the bottom and top bases by $f_1$ and $f_\mathcal{L}$ (without the link) respectively. All $f$, $h_\mathcal{L}$ and $f_\mathcal{L}$ (not counting the link to be built) are Lipschitz-bounded independent on $\mathcal{L}$, hence the added link is as well.

That concludes the construction of the cubical version $f_\mathcal{L}$ of mapping $f_1$. By now there is a potential discrepancy in the definitions for the scale 1, but it is easy to resolve. Just pick the pool of mapping $S_{cell}$ in the definition of the cubical on scale 1 mappings so that they fit into the latter general definition (i.e. their wires are degree-1 cables factoring through the template). We conclude by summarizing the properties of cubical mappings established in this section.

\begin{lemma}
For any cubical on scale 1 mapping $f_1:[0,n\mathcal{L}]^3\rightarrow S^2$ there exists a cubical on scale $\mathcal{L}$ mapping $f_\mathcal{L}$ so that
\begin{enumerate}
\item On each $\mathcal{L}\times \mathcal{L}\times \mathcal{L}$ cell $f_\mathcal{L}$ is given by cables.\\
\item These cables are the ones of a balanced link in the middle of the cell, and the ones factoring through one of the templates. \\
\item Mapping $f_\mathcal{L}$ is homotopic to $f_1$ via a homotopy $h_\mathcal{L}$.\\
\item Homotopy $h_\mathcal{L}$ sends 1-skeleton of $\mathcal{L}$-grid to a point and is ${\sim}1$-Lipschitz over time $\mathcal{L}$ on the rest of the grid.
\end{enumerate} 
\end{lemma} 

One may refer again to the illustration we used in the beginning that we repeat here. It depicts the way the cables can go in one of the templates (but their ends are shrunken to allow to see through the cube). In the bottom right sits a balanced link, moved there also for visibility purposes. 

\includegraphics[scale=0.33]{cableblack.png}

\section{Proof of the main lemma}
\label{proof}

In this section we prove Lemma~\ref{main}. 

First let's overview the idea of the proof. Given a mapping that is cubical on some grid with cells $\mathcal{L}\times \mathcal{L}\times\mathcal{L}$ we want to merge those cells into cells $(2\mathcal{L})^{\times 3}$ and pass to a coarser grid. First we perform a such a merging homotopy on the 2-skeleton (of the coarser grid) rel 1-skeleton. Now on each 3 cell we have some mapping that is homotopic rel boundary to the intended $f_{2\mathcal{L}}$. The difference will be split into several standard Hopf links in a bounded way using the simple large-scale geometry of cubical mappings. We merge all the produced links into a single one, which is trivial due to its Hopf invariant. This nullhomotopy of the difference provides an extension of the homotopy from 2-skeleton to each 3-cell.

\subsection{Boundary homotopy}

The first stage is a homotopy that fixes the mapping on the 2-skeleton of the bigger grid.

Consider a cubical version $f_{2\mathcal{L}}$ of the initial mapping $f_{1}$. Each face $[0,2\mathcal{L}]^2$ has 4 quadrants of the previous grid, in each lives a cable-end-like mapping of $f_{\mathcal{L}}$. We want a homotopy $H_\mathcal{L}:[0,2\mathcal{L}]^2\times [0,T]\rightarrow S^2$ ($T=2\mathcal{L}$) that is constant on the boundary, ends in a bigger cable-end-like mapping $f_{2\mathcal{L}}$ and should be realized by cables. This is a situation similar to the one in definition of cubical mapping: we want to extend cables inside a cube from their ends on the faces. So this is solved in a completely analogous way. We fix appropriate cable templates and put through them cables (using Proposition~\ref{reshape}) so that they cancel portions on the boundary with opposite signs.
This construction is illustrated schematically in the picture below.

\includegraphics[scale=0.2]{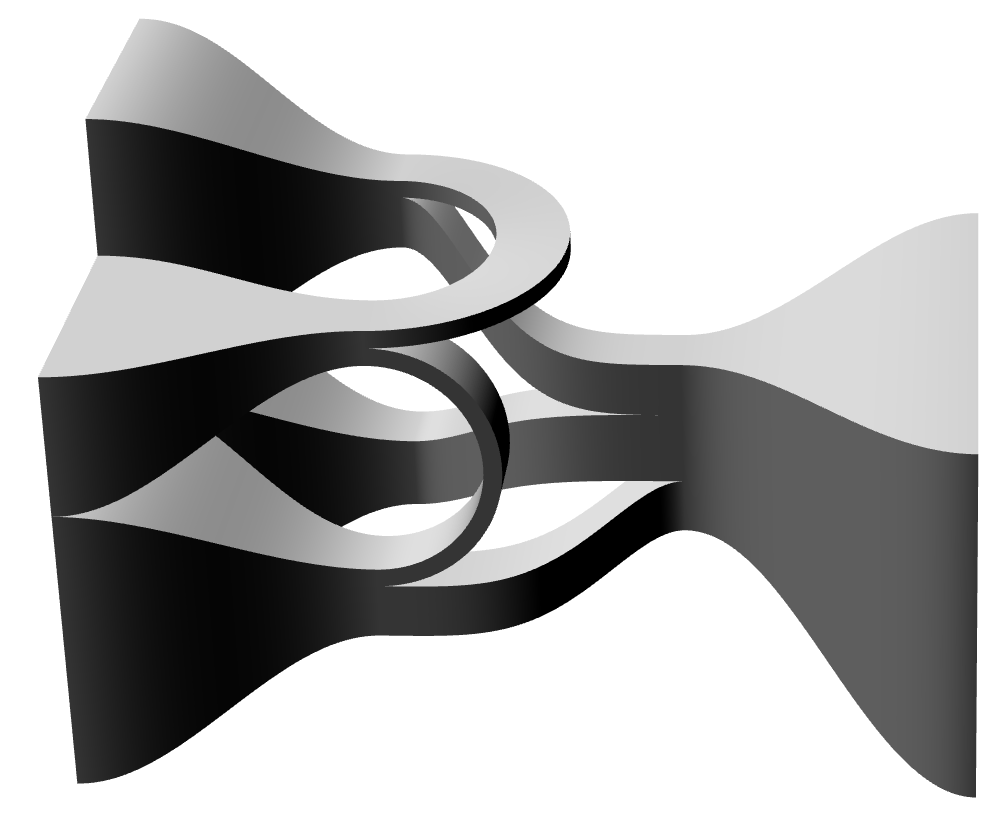}

In order to ensure that the homotopy is extendable we may have to adjust Hopf invariant of $H_\mathcal{L}$. This is done by adding a balanced link $a|b+c|1+\overline{(1)}$ inside the cube $[0,2\mathcal{L}]^2\times [0,T]$. The link should equalize the (relative) Hopf invariant of $H_\mathcal{L}$ to that of some extendable homotopy. We pick the one that is a concatenation of two constructed so far homotopies: one going from $f_\mathcal{L}$ to $f_1$ (that is, all built so far $H_\bullet$ in reverse), and the other being $h_{2\mathcal{L}}$ that goes from $f_1$ to $f_{2\mathcal{L}}$. Since those homotopies are of total size $\lesssim \mathcal{L}$, so is the link that we add to balance out the Hopf invariant.

Extending $H_\mathcal{L}$ to 3-cells is equivalent to contracting the mapping on the boundary of each cell $[0,2\mathcal{L}]^3\times [0,T]$. That is, we consider the clutching mapping that pairs $f_\mathcal{L}$ and $f_{2\mathcal{L}}$ along $H_\mathcal{L}$:
$$f_\mathcal{L}\underset{H_\mathcal{L}}{\sqcup}f_{2\mathcal{L}}:\partial ([0,2\mathcal{L}]^3 \times [0,T])\rightarrow S^2$$ 
that is given by $f_\mathcal{L}$ and $f_{2\mathcal{L}}$ on the bottom and top faces $[0,2\mathcal{L}]^3\times \partial [0,T]$ and by $H_\mathcal{L}$ on the side faces $\partial [0,2\mathcal{L}]^3\times [0,T]$.

\subsection{Null-homotopy of the clutching mapping}

\ \newline

The plan is to treat the few cables (that the mapping is packed into) as individual thickened 1-d objects in $\partial([0,2\mathcal{L}]^3\times [0,T])\cong S^3$, homotope them into few links and then add all these cable links to get a 0 link. It would be that straightforward if only all the cables formed closed non-intersecting loops. But the cables mostly go only from a face to a face and there they split their wires among the cables on the other side of the face. We plan to get rid of this intertwining of wires among cables.

Currently the mapping on $S^3$ can be illustrated by the following picture (which for the sake of readability is of 1 dimension lower, depicting $S^2$ instead, but we will pretend it is $S^3$). 

\includegraphics[scale=0.3]{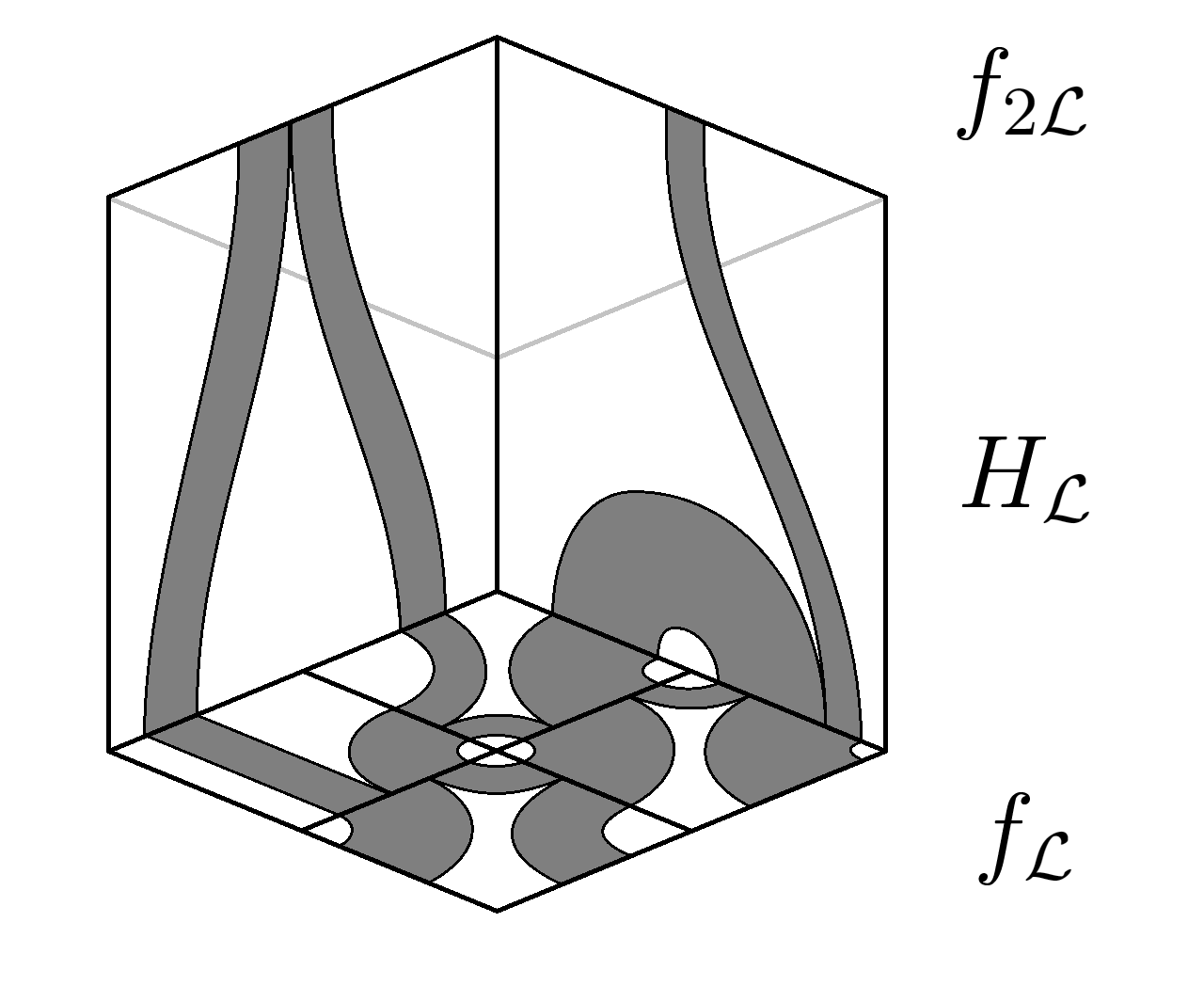}

The lower 3-cell is subdivided in sub-cells of size $\mathcal{L}$ that carry $f_\mathcal{L}$, side cells carry the homotopy $H_\mathcal{L}$ and the top cell of size $2\mathcal{L}$ carries $f_{2\mathcal{L}}$.

A little remark on a sleight of hand that we will use here with our terminology. In the definitions we have required that cables and wires have their ends at the boundary of the domain so that they define continuous mappings to $S^2$. But notice that the condition ``wires end at the boundary'' can be relaxed by adding {``\dots or} the end of the wire matches the start of another wire that continues the first one'' --- that still would generate mappings that are Lipschitz, but will allow us to consider wires (and similarly cables) that end in the interior of the domain. 

Also, rescale the sphere to be of unit size to make the construction uniform. 

The first issue that we deal with, is that right now cables start and end at various 2d-boundaries between 3-cells --- let's call $E_i$ these terminal cross-sections of cable templates (on the previous picture $E_i$ are the intersections of the cable templates with black straight lines). Instead we would like $E_i$ to lie in a single common cross-section. To do so, informally, we just grab the templates by $E_i$ and drag them together, arranging $E_i$ into a single plane, one next the other. Here is a more detailed description.

We start by separating 3-cells apart a bit, as on the next picture, so that gaps of size ${\sim}1$ emerge between cells, and in these gaps the cross-sections $E_i$ get stretched into cylinders (shaded slightly lighter in the picture).

\includegraphics[scale=0.3]{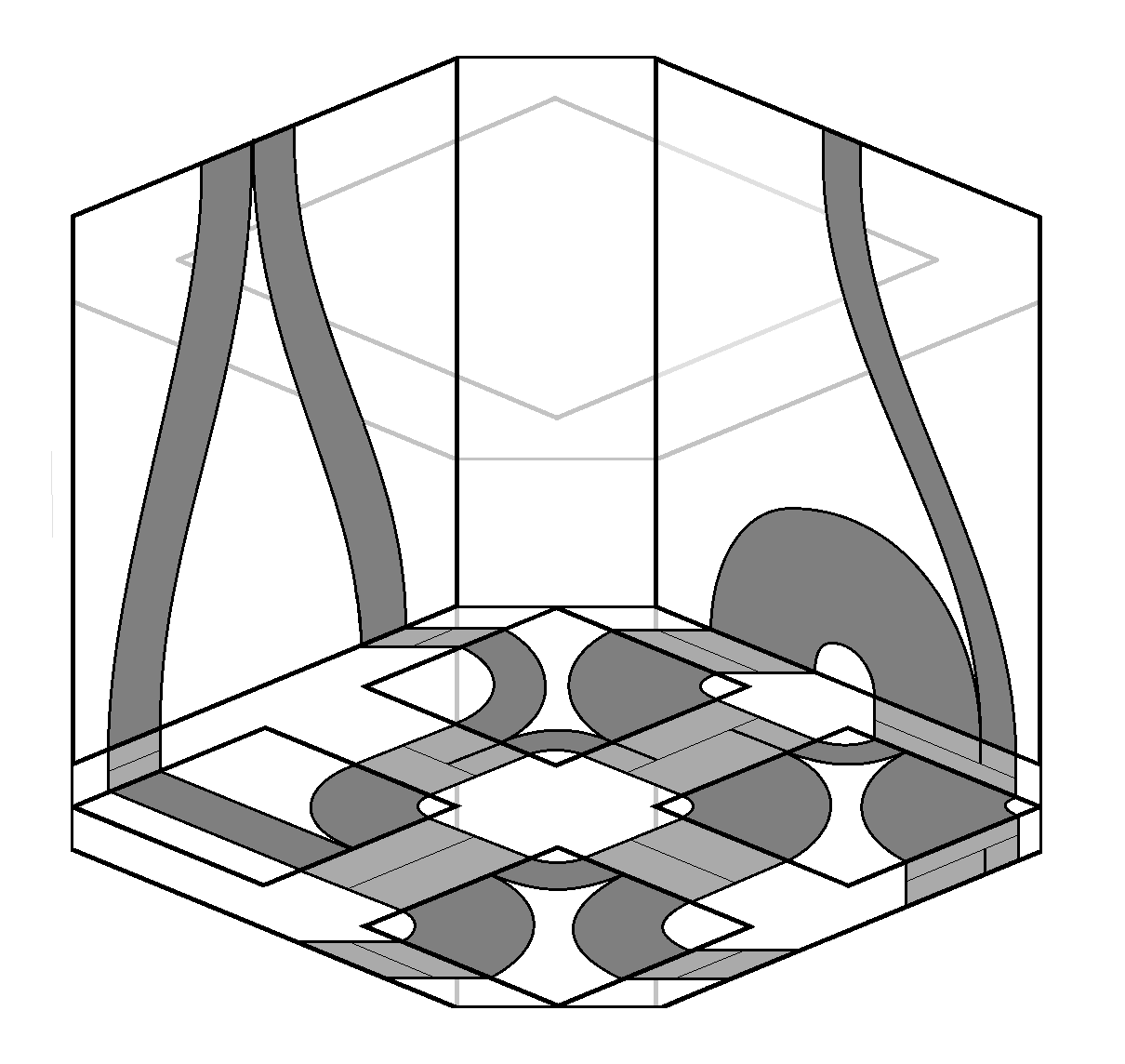}

Pick a ball $B\subset S^3$ of size ${\sim}1$ that is disjoint from the templates (after crating gaps one can fit such a ball in the center of the bottom face). Since each $E_i$ was stretched into a cylinder, let's name the middle section of $E_i$ as $E_i'$ (the thinnest lines on the picture). Pick non-intersecting collection of solid paths $\gamma_i$ of thickness ${\sim}1$ that go (avoiding templates) from the side of $E_i$ to the inside of $B$ so that the ends of $\gamma_i$ in $B$ are placed in line one after another. Now apply a Lipschitz homotopy supported on the union of $\gamma_i$ and $E_i$ that drags $E_i'$ through corresponding $\gamma_i$ till the end, thus placing all the cable terminal section next to each other inside $B$. This process is illustrated in the next picture.
 
\includegraphics[scale=0.3]{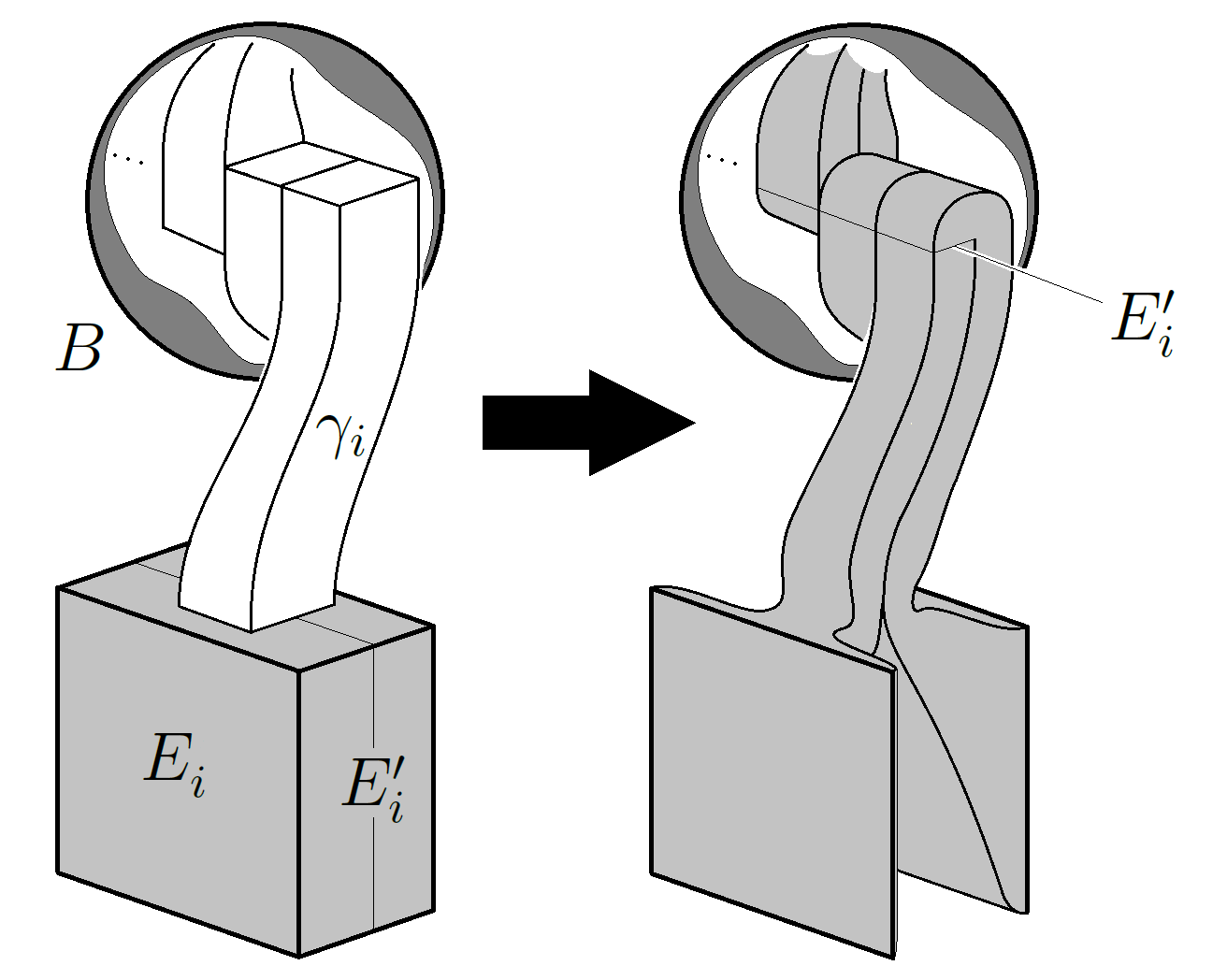}

We arrange the placements of $E_i'$ in agreement with each other. That means they are oriented to have degree of the same sign in that plane, share the coordinates used for ordering the stripes in cable cross-section and are following one another with respect to that order. So that we can imagine the whole collection $E'_i$ as being one ``meta-cable'' cross-section, whose wires then split into smaller cables, $i$-th cable being composed of consecutive wires stating from $c_i$-th up to $c_{i+1}$-th. Those cables leave this meta-cable cross-section $\bigcup E'_i$, wander around in $S^3$ and end now at the same section $\bigcup E'_i$.  

The second issue to resolve now is that once the $i$-th cable comes back, it reconnects its wires $[c_i,c_{i+1})$ not necessarily to themselves but to some $[c_i+n_i,c_{n+1}+n_i)$ (preserving their respective order, however). We would like this monodromy of tracking wires along the cables to be trivial, so that the meta-cable would split into disjoint circular cables. To trivialize the monodromy in a bounded manner we use the following proposition.

\begin{proposition}
Let a permutation $\sigma\in S_N$ be such that it rearranges $k$ blocks. That is, $[1,N]$ is a disjoint union of $k$ blocks $[c_i,c_{i+1})$ and there are $n_i\in\mathbb{Z}$ such that for $x\in [c_i,c_{i+1})$ one has 
$$\sigma(x)=x+n_i.$$
Then $\sigma$ is a composition of $k-1$ permutations $\sigma_i$ each of which rearranges just 2 blocks. 
\end{proposition}

\begin{proof}
It follows by induction. Consider the composition $\sigma'$ of $\sigma$ and the swapping that shifts ${[\sigma(1),N]}$ to ${[1,N-(\sigma(1)-1)]}$ and shifts ${[1,\sigma(1))}$ to ${[N-\sigma(1)+2,N]}$. The permutation $\sigma'$ leaves intact the first block $[1,c_2)$ and rearranges the rest $k-1$ of them. So it is by assumption a composition of $k-2$ swappings. That gives a decomposition of $\sigma$ into $k-1$ swappings.
\end{proof}

Applying the proposition to the monodromy of the wires in meta-cable, we represent it as a bounded composition of transpositions of the type ``split wires into 2 blocks and swap them". Each such operation is realizable by a {\bf swapping cable} of degree $d$ that goes in a loop but midway splits in 2 sub-cables of degrees $d_1$ and $d_2$ that swap their places and reunite, as in the next figure.

\includegraphics[scale=0.3]{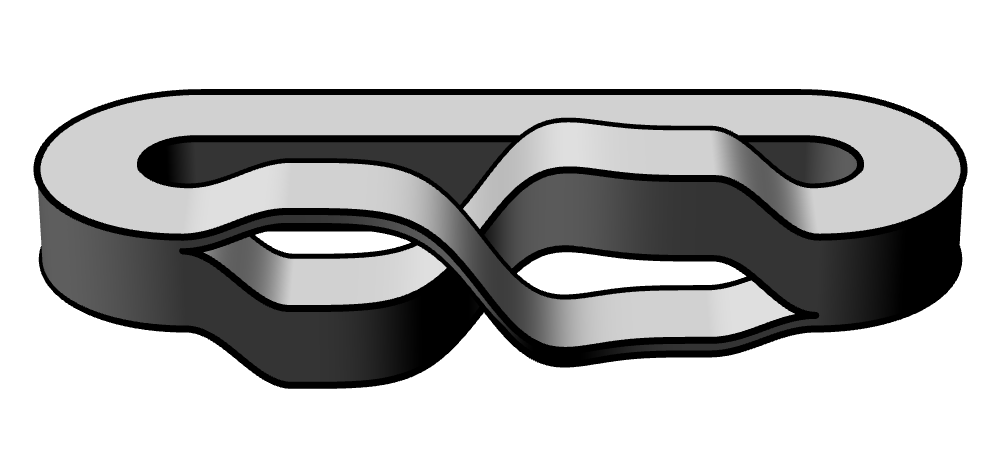}

So, by adding the appropriate swapping cables to a constant piece of the meta-cable we trivialize the monodromy in wires. To add a swapping cable we first create a pair of it and its opposite and then insert one of them, leaving the other outside. 

The byproduct swapping cables themselves can be broken down to pieces we are familiar with. Decompose the swapping of sub-cables into 1) revolving their pair as a solid and 2) revolving them back separately on each own to compensate the twist gained in 1). This breaks the swapping cable into a combination of 3 twisted links $\underline{d}-(\underline{d_1}+\underline{d_2})$ as shown in the next figures.

\includegraphics[scale=0.3]{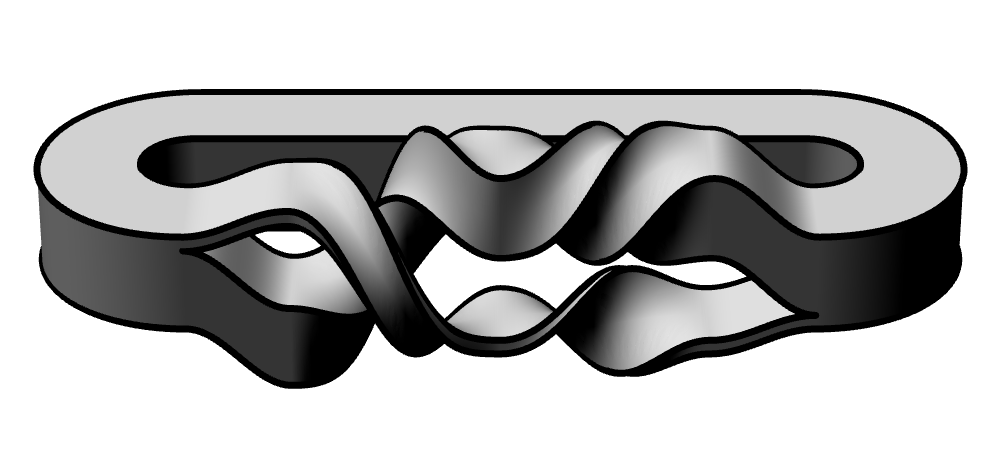}

\includegraphics[scale=0.3]{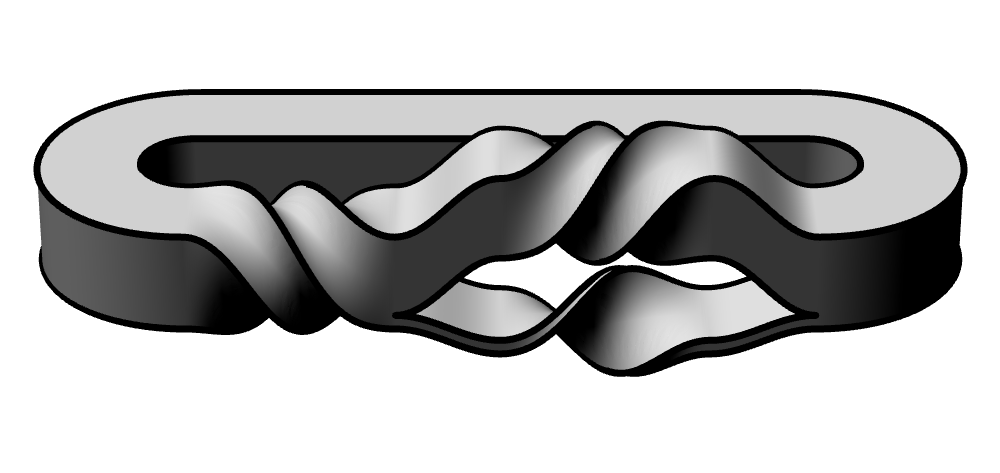}

Each of the swapping cables is broken this way into three twisted links, which are further broken into balanced links by Corollary~\ref{twist-bal}. 

Once the monodromy of the meta-cable is trivial, it can be viewed just as a collection of separate tangled cables it was built from, that close up now perfectly to themselves with no permutations. With aid of propositions~\ref{comb} and~\ref{reshape} we can even assume that they have constant cross-sections we favor. Namely, Proposition~\ref{reshape} (with cable length regarded as time) allows to put one cross-section of a cable into a desired shape and then Proposition~\ref{comb} allows to homotope the whole cable to have the same cross-section. 

Now, contract these cables one by one to the state of a localized loop $d_i|0$, maybe producing on the way 1) a bounded number of links $d_i|d_j$ for every time $i$-th cable has to pass across $j$-th cable, and 2) several twisted links to deal with the twisting that the cable might have in the end (one full turn of a cable is two twists $\underline{d_i}$ plus a turn of each individual wire $(d_i\times\overline{1})\backsimeq ([d_i/2]|1+\overline{(1)})\backsimeq$ a balanced link). The amount of links produced is bounded, since the arrangements of the cables came from a finite set of template choices for the initial 8 sub-cells, the final bigger set and the homotopies $H'_\mathcal{L}$, and for each arrangement the constructed contractions and untwistings produce a finite number of new links. 

After each cable is contracted, we are left with a bounded number of balanced links of size ${\lesssim}\mathcal{L}$ that were produced at various stages of the homotopy. So we can add them up into a single link by Lemma~\ref{ad}, which happens to be trivial due to vanishing Hopf invariant. Thus we constructed a bounded null-homotopy of the difference mapping and thus finish establishing a bounded homotopy the Lemma~\ref{main}.

\section{Proof of the main theorem}

Here we show in detail how Lemma~\ref{main} implies Theorem~\ref{mainth}.

First we switch our domain from the whole $S^3$ to a cube $[0,1]^3$ that is ${\sim 1}$-biLipschitz embedded into $S^3$. This is validated by means of precomposition with a ${\sim}1$-Lipschitz contraction $H_c$ of the complement to the cube in $S^3$ to the base point. Then we rescale the cube to be $[0,2^N]^3$ where $2^N{\sim}(L+1)$, so that we consider ${\sim}1$-Lipschitz mappings and aim for a homotopy of length ${\sim}(L+1)$.

In the definition of cubical on scale 1 mapping we established that for some $c>0$ all $c$-Lipschitz mappings from a unit grid in $\mathbb{R}^3$ to $S^2$ are ${\sim}1$-Lipschitz homotopic (over unit time) to a cubical mapping. So we pick an integer $N{\sim} \log_2(L+1)$ such that the initial contraction $H_c$ and rescaling of the unit cube to the size $2^N$ multiply the Lipschitz constant of mappings by a factor $<c/L$. Thus starting with $L$-Lipschitz mapping on $S^3$ we get a ${<}c$-Lipschitz mapping on a cube, and that mapping is homotoped to a cubical one. 

Now we use the Lemma~\ref{main} to pass to higher scales. The homotopy that the lemma produces at step $i$ is happening at scale $2^i$ and hence has length ${\sim}2^i$. So catenating all of them gives a homotopy of length $${\sim}(2^0+2^1+\cdots +2^N)=2\cdot 2^N -1\lesssim(L+1)$$
through mappings of Lipschitz constant ${\sim}1$. At this point both end mappings $f_0$ and $f_1$ are reduced to a mappings cubical on a full scale $2^N$. Since the boundary of the $2^N$-cube was sent to the base point at all times, the only cables present in the cubical mappings are those of the balanced links inside the cube. Since the mappings were homotopic by assumption, Proposition~\ref{cancel} implies then existence of a homotopy between the remaining balanced links over time ${\sim}2^N$ through Lip${\sim}1$ mappings.

Patching all these homotopies together we get one of length ${\sim} (L+1)$ through mappings that are ${\sim}(L+1)$-Lipschitz (that is, with regard to the metric of $S^3$ itself). 

\section*{Appendix}
\label{app}
Here we show that if $m\geqslant n$, there exist ${\sim}L$-Lipschitz nullhomotopic mappings $f_L:S^m\rightarrow S^n$ such that the length of minimal null-homotopy of $f_L$ through ${\sim}L$-Lipschitz mappings is at least ${\sim} L$. Indeed, let $f_L$ be  the Whitehead product of a constant mapping and ${\sim}L$-Lipschitz mapping $g_L:S^n\rightarrow S^n$ of degree $L^n$. Indeed, there is a subset ${(D^n\times S^{m-n})}\subset S^m$ so that $f_L$ restricts to it as
$$f_L:\ \ (D^n\times S^{m-n})\stackrel{pr}{\rightarrow} D^n\stackrel{/\partial}{\rightarrow}S^n\stackrel{g_L}{\rightarrow}S^n.$$
Let $h_L:S^m\times[0,1]$ be a null-homotopy of $f_L$. Since the mapping  $h_L$ extends from the boundary $\partial (D^n\times * \times[0,1])$ to the interior, it should have zero degree on the boundary. Its restriction to the bottom $D^n\times * \times \{0\}$ is $g_L$ that is carrying degree $L^n$, restriction to the top is constant, thus carrying degree 0, so it should also have degree of size $L^n$ on the sides $S^{n-1}\times * \times[0,1])$. A homotopy of length $l$ through ${\sim}L$ mappings on this space has degree at most ${\sim}L^{n-1}l$, therefore the length $l$ is at least ${\sim}L$.

\bibliography{ref}
\bibliographystyle{abbrv}
\end{document}